\documentclass[a4paper,11pt]{amsart}
\usepackage{amsmath,amsfonts,amssymb,graphicx,url,amsthm,hyperref}
\usepackage[left=2.5cm,right=2.5cm,top=3cm,bottom=3cm]{geometry}
\usepackage[nocompress]{cite}
\usepackage{xcolor}
\usepackage{algorithm}
\usepackage[noend]{algpseudocode}
\algrenewcommand\algorithmicrequire{\textbf{Input:}}
\algrenewcommand\algorithmicensure{\textbf{Output:}}
\usepackage[T1]{fontenc}
\usepackage{comment}
\usepackage{enumerate}
\usepackage{enumitem}
\usepackage{graphicx}
\usepackage{subcaption}
\usepackage{booktabs}
\usepackage{multirow}
\usepackage{tikz}
\usetikzlibrary{arrows.meta, positioning}

\newcommand{\kernel}{\mathcal{K}}

\newcommand{\bs}{\boldsymbol}
\newcommand{\bm}{\bs}

\newcommand\Fnor[1]{F^{\lambda}_{\mathrm{nor}}(#1)}
\newcommand\Fnat[1]{F^{\lambda}(#1)}
\newcommand{\Act}{\mathcal{A}}
\newcommand{\Inact}{\mathcal{I}}
\newcommand{\mer}{H_\tau}
\newcommand{\iprod}[2]{\langle #1, #2 \rangle}

\DeclareMathOperator{\diam}{diam}

\DeclareMathOperator{\spn}{span}

\DeclareMathOperator\dist{dist}
\DeclareMathOperator\supp{supp}

\usepackage{tikz,pgfplots,pgfplotstable}
\pgfplotsset{compat=newest}
\usetikzlibrary{arrows,decorations.pathmorphing,backgrounds,
positioning,calc,matrix}
\usepackage{graphicx,color,nicefrac} 
\usepackage{amssymb, amsmath, amsbsy,mathtools}
\usepackage{amsmath,amsfonts}
\usepackage[export]{adjustbox}
\newtheorem{theorem}{Theorem}[section]

\newtheorem{remark}[theorem]{Remark}

\def\letters{a,b,c,d,e,f,g,h,i,j,k,l,m,n,o,p,q,r,s,t,u,v,w,x,y,z}
\def\Letters{A,B,C,D,E,F,G,H,I,J,K,L,M,N,O,P,Q,R,S,T,U,V,W,X,Y,Z}
\makeatletter
\@for \@l:=\Letters \do{%
  \expandafter\edef\csname\@l bb\endcsname{\noexpand\ensuremath{%
  \noexpand\mathbb{\@l}}}%
  \expandafter\edef\csname\@l bf\endcsname{{\noexpand\bf \@l}}%
  \expandafter\edef\csname\@l cal\endcsname{\noexpand\ensuremath{%
  \noexpand\mathcal{\@l}}}%
  \expandafter\edef\csname\@l eu\endcsname{\noexpand\ensuremath{%
  \noexpand\EuScript{\@l}}}%
  \expandafter\edef\csname\@l frak\endcsname{\noexpand\ensuremath{%
  \noexpand\mathfrak{\@l}}}%
  \expandafter\edef\csname\@l rm\endcsname{{\noexpand\rm \@l}}%
  \expandafter\edef\csname\@l scr\endcsname{\noexpand\ensuremath{%
  \noexpand\mathscr{\@l}}}%
}
\@for \@l:=\letters \do{%
  \expandafter\edef\csname\@l bf\endcsname{{\noexpand\bf \@l}}%
  \expandafter\edef\csname\@l frak\endcsname{\noexpand\ensuremath{%
  \noexpand\mathfrak{\@l}}}%
  \expandafter\edef\csname\@l scr\endcsname{\noexpand\ensuremath{%
  \noexpand\mathscr{\@l}}}%
}
\makeatother
\definecolor{shadecolor}{rgb}{0.6, 0.6, 0.6} 
\definecolor{red}{rgb}{1,0,0} 
\definecolor{darkgreen}{rgb}{0, 0.6, 0}
\definecolor{darkred}{rgb}{0.9,0,0} 
\newcommand{\isdef}{\mathrel{\mathrel{\mathop:}=}}
\newcommand{\defis}{\mathrel{=\mathrel{\mathop:}}}

\begin{document}
\title[Tree-Adaptive Multiscale Kernel Lasso 
in Samplet Coordinates]{%
Tree-Adaptive Multiscale Kernel Lasso 
in Samplet Coordinates}
\author[
S. Avesani,
G. Fumagalli,
M. Multerer,
and C. Segala
]{
Sara Avesani$^*$,
	Gaia Fumagalli$^*$,
	Michael Multerer$^*$,
	and Chiara Segala$^*$
}
\thanks{$^*$IDSIA USI-SUPSI,
  Universit{\`a} della Svizzera italiana,
	Via la Santa 1, 6962 Lugano, Svizzera.\\
\{\texttt{sara.avesani, gaia.fumagalli, michael.multerer,
chiara.segala}\}\texttt{@usi.ch}}

\keywords{kernel-based approximation, multiresolution analysis, samplets, adaptive
	subsampling, basis pursuit, sparse regularization}

\subjclass[2020]{46E22, 65D15, 65K10, 41A30}

\begin{abstract}
We develop a novel framework for sparse multiscale kernel approximation of large 
scattered data problems based on a samplet representation. Samplets form a
multiresolution analysis of localized discrete signed measures and enable
quasi-sparse representations of kernel matrices associated to asymptotically
smooth kernels as well as smoothness detection of scattered data.
Building on the latter, we introduce an adaptive data site selection strategy
based on the native reproducing kernel Hilbert space norm of the samplet
expansion coefficients.
The selection results in a small set of
representative data sites, significantly reducing the effective problem size. 
On the corresponding reduced kernel subspace, we solve an $\ell^1$-regularized
least-squares problem using a trust-region semismooth Newton method in a
normal-map formulation, stabilized by an online low-rank SVD on the active set
to handle the notorious ill-conditioning of kernel matrices. Numerical
experiments in two and three dimensions, including multi-kernel models with
varying lengthscales, demonstrate that the proposed approach achieves accurate
reconstructions with considerably sparser representations and good computational
efficiency.
\end{abstract}

\maketitle



\section{Introduction}\label{sec:intro}
Kernel-based methods are a powerful tool for the approximation of
functions from scattered data \cite{Wendland2004,Fasshauer2007}, with
applications ranging from uncertainty quantification 
\cite{HMQ25,jakeman2015enhancing} and PDE-constrained optimization 
\cite{schaback2006kernel}
to financial econometrics and machine learning 
\cite{scholkopf2002learning,Rasmussen2006}.
Their feasibility for large-scale problems is, however, limited by the dense
and ill-conditioned nature of kernel matrices, which leads to prohibitive
computational cost and memory requirements.

The computational cost as well as memory requirements can be mitigated by
compression techniques for kernel matrices. Existing approaches comprise 
low-rank techniques \cite{WS00, RR17,FMS25} as well as 
block low-rank techniques such as the fast multipole method
\cite{greengard1987fast} or $\mathcal{H}/\mathcal{H}^2$-matrices
\cite{Hac15,borm2010efficient}. For asymptotically smooth kernels,
samplets, which form a multiresolution analysis of discrete signed measures,
see \cite{HM22}, yield quasi-sparse representations of the corresponding
kernel matrices that can be compressed to essentially linearly many remaining 
entries in terms of the number of data sites.
In contrast to classical wavelet or Fourier constructions, see, for example,
\cite{Mallat}, samplets are specifically tailored
to scattered data, making them particularly effective for
kernel-based approximation.

A second, complementary mechanism exploits the multiscale organization of the
data to reduce the number of degrees of freedom.
Within the samplet structure, which is organized in terms of a cluster tree,
the information content of the data is
naturally localized in the clusters at the different scales.
We leverage this structure to introduce a tree-adaptive subsampling strategy
that selects only a small set of representative data sites prior to solving
the kernel regression problem.
The relevance of each cluster can be quantified through {the kernel's native
reproducing kernel Hilbert space norm of its samplet coefficients}.
This yields a geometrically meaningful, kernel-dependent measure of importance,
which accounts for correlation and smoothness effects.
Following the adaptive tree selection strategy of \cite{binev2004fast},
cluster contributions are aggregated bottom-up and propagated top-down,
leading to the identification of a high-energy subtree. Selecting representative
data sites from each of the remaining clusters, the resulting reduced set
reflects the intrinsic multiscale complexity of the data, rather than its raw
size, and enables a substantial reduction of the effective problem dimension.
{The tree approximation algorithm of
\cite{binev2004fast} computes a near-optimal adaptive subtree in
linear cost, where the energy of a cluster is the squared Euclidean norm
of the coefficients supported on it, so that the resulting subtree is
tailored to the $\ell^2$-norm of the coefficient vector. 
We propose a kernel-weighted energy that accounts
for the correlations between data sites induced by the kernel {and bounds
the native space norm from above}. The Euclidean energy of \cite{binev2004fast}
is recovered as a special case.}

On the adaptively reduced subspace of kernel translates
corresponding to the selected data sites, we consider
$\ell^1$-regularized least-squares formulations to promote sparse
representations of the solution.
Such formulations can be interpreted in the spirit of basis pursuit and
sparse approximation in kernel dictionaries
\cite{parikh2014proximal,smola2000sparse,tropp2004greed},
and have recently been studied in the samplet setting
\cite{BaroliHarbrechtMulterer2024}.
Beyond sparsity, the $\ell^1$-penalty provides a mechanism to automatically
select relevant kernel contributions, which is particularly effective in
multi-kernel settings combining kernels with varying lengthscales.
{The resulting optimization problem is solved by a trust-region semismooth
Newton method \cite{qi1993nonsmooth,HIK,Ulbrich2001} based on a normal-map
formulation \cite{ouyang2023trust,ouyang2025trust}. This globalized
second-order method for composite problems applies semismooth Newton steps
to the normal map and accepts or rejects them through a trust-region
mechanism. In the $\ell^1$-regularized setting, the proximity operator is
the soft-shrinkage operator, so that the Newton system reduces to the
active set. This reduced system is solved by a block online low-rank SVD update
of the active kernel block, see \cite{brand2002incremental}.}

The main objective of this work is to develop a multiscale framework for
sparse kernel approximation that combines samplet-based compression with
kernel-driven adaptive subsampling. By first reducing the set of data sites
and then promoting sparsity in the
coefficients, the proposed approach achieves accurate multiscale
reconstructions with significantly reduced computational cost.
The combination of samplet compression, adaptive data site selection
and sparse regression enables efficient treatment of large-scale,
heterogeneous data sets.

The remainder of this paper is organized as follows.
Section~\ref{section:Problem_formulation} introduces the kernel approximation
framework and the samplet representation, together with the adaptive
subsampling strategy. Section~\ref{sec:newton_method} recalls the
$\ell^1$-regularized formulation together with the proposed solution algorithm.
Section~\ref{sec:results} reports numerical experiments validating
the proposed approach on problems of increasing complexity. Concluding remarks
are stated in Section~\ref{sct:conclusion}.

\section{Kernel-based approximation framework}
\label{section:Problem_formulation}

In this section, we introduce the kernel-based approximation framework
underlying the proposed approach. We briefly recall the main concepts of
reproducing kernel Hilbert spaces, describe the samplet-based
compression of kernel matrices, and present the tree-adaptive
subsampling strategy used to construct reduced approximation subspaces.

\subsection{Reproducing kernel Hilbert spaces}
Let \((\Hcal,\langle\cdot,\cdot\rangle_\Hcal)\) be a Hilbert space of
continuous, real-valued functions \(h\colon\Omega\to\Rbb\) defined on a bounded
Lipschitz region \(\Omega\subset\Rbb^d\). Then, 
\((\Hcal,\langle\cdot,\cdot\rangle_\Hcal)\) is a
\emph{reproducing kernel Hilbert space} and exhibits a
\emph{reproducing kernel} $\kernel \colon \Omega \times \Omega \to \Rbb$
with \(\kernel({\bm x},\cdot)\in\Hcal\) 
for every \({\bm x}\in\Omega\) and \(h({\bm x})=\langle\kernel
({\bm x},\cdot),h\rangle_\Hcal\) for every \(h\in\Hcal\). 
The reproducing kernel is symmetric and positive definite in the sense
that the associated 
\emph{kernel matrix} 
\begin{equation}\label{eq:KernelMatrix}
{\bm K}\isdef [\kernel({\bm x}_i,{\bm x}_j)]_{i,j=1}^N\in\Rbb^{N\times N}
\end{equation}
is symmetric and positive definite for any set of mutually distinct points 
\({\bm x}_1, \ldots,{\bm x}_N\in\Omega\) and any \(N\in\Nbb\).

For \(h\in\Hcal\) and a $\bm x \in \Omega$, we denote the point evaluation
functional by
\(
\delta_{{\bm x}}(h)\isdef h({\bm x}),\ {\bm x}\in\Omega
\)
and notice that there holds $\delta_{\bm x}\in\Hcal'$ for any
${\bm x}\in\Omega$. By the reproducing property, the
kernel \(\kernel({\bm x},\cdot)\) is then nothing but the Riesz representer
of $\delta_{\bm x}$, 
i.e.,
\(
\delta_{{\bm x}}(h)=\langle\kernel({\bm x},\cdot),h\rangle_\Hcal
\ \text{for every } h\in\Hcal.
\)
Given the \emph{data sites}
\(
X=\{{\bm x}_1,\ldots,{\bm x}_N\}\subset\Omega,
\)
we define the finite-dimensional subspace
\begin{equation}\label{eq:HX}
\Hcal_X \isdef \spn\{\phi_1,\ldots,\phi_N\} \subset\Hcal,
\end{equation}
 of \emph{kernel translates}
 $\phi_i \isdef \kernel({\bm x}_i,\cdot)$,
$i=1,\ldots,N$, associated to $X$. By the Riesz isometry,
the subspace \(\Hcal_X\) is isometrically
isomorphic to the subspace 
\(
    \Hcal_X'\isdef\spn\{\delta_{{\bm x}_1},\ldots,
\delta_{{\bm x}_N}\}\subset\Hcal'.
\)
More precisely, for any coefficient vector $\bs\alpha =
[\alpha_i]_{i=1}^N \in \Rbb^N$,
we identify
\[
h'=\sum_{i=1}^N\alpha_i\delta_{{\bm x}_i}\in\Hcal_X'\quad\text{and}\quad
h=\sum_{i=1}^N\alpha_i\kernel({\bm x}_i,\cdot)\in\Hcal_X.
\]
Especially, there holds
\(
  \langle \delta_{{\bm x}_i},\delta_{{\bm x}_j}\rangle_{\Hcal'}
  =
    \langle \kernel({\bm x}_i,\cdot),\kernel({\bm x}_j,\cdot)\rangle_{\Hcal}
    =\kernel({\bm x}_i,{\bm x}_j),
\)
such that the Gramian of either basis becomes the kernel matrix
\eqref{eq:KernelMatrix}. 
The duality between \(\Hcal_X\)
and  \(\Hcal_X'\) implies that
the \(\Hcal\)-orthogonal projection
of a function \(h\in\Hcal\) onto \(\Hcal_X\) is given by the kernel interpolant
\begin{equation*}
s_h\isdef\sum_{i=1}^N\alpha_i\kernel({\bm x}_i,\cdot),
\end{equation*}
which satisfies the interpolation conditions
\(
s_h({\bm x}_i) = h({\bm x}_i)
\)
for $i=1,\ldots, N$.
Equivalently, $s_h$ is characterized by the variational formulation
\begin{equation}\label{eq:Galerkin}
\langle s_h,v\rangle_\Hcal
=\langle h,v\rangle_\Hcal\quad\text{for all }v\in\Hcal_X.
\end{equation}
This is seen by choosing the basis of kernel translates 
as ansatz and test functions in \eqref{eq:Galerkin}.
The expansion coefficients 
\({\bm\alpha}\isdef[\alpha_i]_{i=1}^N
\)
can be retrieved by solving the linear system
\begin{equation}\label{eq:LSE}
{\bm K}{\bm\alpha}={\bm h},\quad{\bm h}\isdef[h({\bm x}_i)]_{i=1}^N.
\end{equation}
For a growing number of data sites, assembly and solution of \eqref{eq:LSE}
easily become prohibitive. To speed up computations, we transform the
system using a multiresolution basis resulting in a quasi-sparse
kernel matrix.

\subsection{Samplets and kernel matrix compression}\label{sec:samplets}

Samplets are discrete, localized signed measures that exhibit vanishing
moments. We briefly recall their underlying concepts as introduced in
\cite{HM22}. To this end, we shall consider the space 
\(\Xcal'\isdef\Hcal_X'\) equipped with a different inner product
\(\langle\cdot,\cdot\rangle_{\Xcal'}\) such
that \(\langle\delta_{{\bs x}_i},\delta_{{\bs x}_j}
\rangle_{\Xcal'}=\delta_{i,j}\), see \cite{BM24} for details.
Then, there holds
\begin{equation*}
\langle u',v'\rangle_{\Xcal'}\isdef\sum_{i=1}^N u_iv_i,\quad
u'=\sum_{i=1}^Nu_i\delta_{{\bm x}_i},\ v'=\sum_{i=1}^Nv_i\delta_{{\bm x}_i},
\end{equation*}
and we retrieve the Euclidean geometry for the coefficient vectors.
This is in contrast to the restriction of the Hilbert space inner product
$\langle \cdot,\cdot\rangle_\Hcal$ to \(\Hcal_X\), which is rather expressed in
terms of the kernel matrix $\bm K$ from \eqref{eq:KernelMatrix} than by the
identity matrix. Specifically, for $ u, v \in \Hcal_X$ represented in the 
basis of kernel translates with coefficient vectors
\({\bm u}\isdef[u_i]_{i=1}^N\) and \({\bm v}\isdef[v_i]_{i=1}^N\), we have
\(
\langle u,v\rangle_\Hcal={\bm u}^\intercal{\bm K}{\bm v}. 
\)
In this regard, the spaces \((\Xcal',\langle\cdot,\cdot\rangle_{\Xcal'})\)
and \((\Hcal_X',\langle\cdot,\cdot\rangle_{\Hcal'})\) are equivalent, where
the norm equivalence constants are the smallest
and the largest eigenvalue of the kernel matrix, respectively.

The construction of samplets is based on a multiresolution analysis of the space
\(\Xcal'\), i.e., a nested sequence of subspaces
\(
\Xcal_0'\subset\Xcal_1'\subset\cdots\subset\Xcal_J'=\Xcal',
\)
where \(J\sim\log N\). {Geometrically, this multiresolution analysis is
induced by a \emph{hierarchical cluster tree} \(\Tcal\) for
\(X=\{{\bs x}_1,\ldots,{\bs x}_N\}\), that is, a binary tree with root \(X\)
whose nodes \(\tau\in\Tcal\), called \emph{clusters}, satisfy
\(\tau=\tau' \cup \tau''\) whenever
\(\operatorname{child}(\tau)=\{\tau',\tau''\}\). The \emph{level} \(j\) of a
cluster is its distance to the root, while clusters without children are the
\emph{leaves} of \(\Tcal\). Here, \(\Tcal\) is built by recursively splitting
each cluster into two subsets of essentially equal cardinality. We refer to \cite{HM22} for the details.

The multiresolution analysis is then obtained by performing, along the branches
of \(\Tcal\), a suitable change of basis of the point evaluation functionals
associated to each cluster.}
This change of basis can be constructed such that we obtain an
\(\Xcal'\)-orthogonal splitting \(\Xcal_{j+1}'=\Xcal_j'\oplus\Scal_{j+1}'\)
in terms of the \emph{detail spaces}
\(\Scal_{j+1}'\perp\Xcal_j'\). Consequently, there holds
\(\Xcal_J'=\bigoplus_{j=0}^{J}\Scal_{j}'\) with \(\Scal_{0}'\isdef\Xcal_0'\). 
The union of the bases \(\{\sigma_{j,k}\}_k\) of the spaces \(\Scal_j'\) for
\(j=0,\ldots,J\) is called a \emph{samplet basis} for \(\Xcal_J'\).

The following theorem collects the main properties of the samplet basis,
see \cite{HM22, avesani2025multiresolution} for details. 
\begin{theorem}\label{theo:waveletProperties}
The samplet basis \(\bigcup_{j=0}^{J}\{\sigma_{j,k}\}_k\)
forms an orthonormal basis in $\Xcal'$, satisfying
the following properties:
\begin{enumerate}[label=(\roman*)]
\item 
There holds $\operatorname{dim}\Scal'_j\sim 2^{dj}$.
\item 
The samplet basis exhibits \emph{vanishing moments} of order \(q+1\), i.e.,
\begin{equation}\label{eq:vanishingMoments}
\sigma_{j,k}(p)
 = 0\quad\text{for all}\ p\in\Pcal_q,\ j\geq 1,
\end{equation}
where \(\Pcal_q\isdef\operatorname{span}
\{{\bs x}^{\bs\alpha}:\|\bs\alpha\|_1\leq q\}\)
is the space of polynomials of degree at most \(q\).
\item 
The coefficient vector 
${\bm\omega}_{j,k}=\big[\omega_{j,k,i}\big]_i$ 
of the samplet $\sigma_{j,k}$ satisfies 
\(\|{\bm\omega}_{j,k}\|_{1}\lesssim 2^{(J-j)d/2}\).
\item {Let $x_0 \in \Omega$, $f \in C^\alpha(x_0)$ for some $\alpha \geq 0$
with $q \geq \lfloor \alpha \rfloor$. Then, for every cluster 
$\tau_{j,k}\in\Tcal$ containing $x_0$, there holds
\begin{equation*}
|\sigma_{j,k}(f)| \lesssim \operatorname{diam}(\tau_{j,k})^\alpha
\sqrt{\#\tau_{j,k}},
\end{equation*}}
where $C^\alpha(x_0)$ denotes the Jaffard regularity at $x_0$. 
\end{enumerate}
\end{theorem}

Due to the orthonormality of the samplet basis, the samplet transform
\(
[\sigma_{j,k}]_{j,k}={\bs T}[\delta_{{\bs x}_i}]_{i=1}^{N}
\)
has the property \({\bs T}^\intercal{\bs T}={\bs T}{\bs T}^\intercal
={\bs I}\in\Rbb^{N\times N}\). If the samplet transform 
\({\bs T}\) is computed recursively with respect to the hierarchy induced
by \(\Tcal\), its cost is of order \(\Ocal(N)\).

Applying the Riesz isometry to \(\{\sigma_{j,k}\}_{j,k}\) gives rise to
a basis in \(\Hcal_X\) defined as
\[
\psi_{j,k}\isdef\sum_{i=1}^N\omega_{j,k,i}\,\kernel({\bm x}_i,\cdot)\in\Hcal_X.
\]
The Gramian of this transformed samplet basis is obviously given by
\begin{equation}\label{eq:K-samplet-basis}
\big[\langle\psi_{j,k},\psi_{j',k'}\rangle_\Hcal\big]_{j,k,j',k'}
= {\bm T}{\bm K}{\bm T}^\intercal
\defis {\bm K}^\Sigma.
\end{equation}
Therefore, solving the linear system \({\bm K}{\bm\alpha}={\bm h}\) arising from
\eqref{eq:Galerkin}, is equivalent to solving
\begin{equation} \label{eq:system_samplets}
{\bm K}^\Sigma{\bm\beta} = {\bm T}{\bm h},
\qquad 
{\bm\beta} = {\bm T}{\bm\alpha},
\end{equation}
Furthermore, the basis \(\{\psi_{j,k}\}_{j,k}\) inherits the vanishing
moment property. Specifically, for any \(h\in\Hcal\) such that
\(h|_{\supp\sigma_{j,k}}=p|_{\supp\sigma_{j,k}}\) for some \(p\in\Pcal_q\),
the vanishing moment property \eqref{eq:vanishingMoments} implies
\(
\langle\psi_{j,k},h\rangle_\Hcal = 0.
\)
The condition \(h|_{\supp\sigma_{j,k}}=p|_{\supp\sigma_{j,k}}\) can be
translated into the requirement that there exists a polynomial of
total degree \(q\) interpolating all points in the support of \(\sigma_{j,k}\).
Thus, if the coefficient vector of \(h\) in the basis of kernel translates
of \(\Hcal_X\) can locally be interpolated by a polynomial of degree at
most \(q\), the corresponding samplet coefficients are zero. For such
functions \(h\), we obtain sparse representations in the samplet basis.

In what follows, we assume that $\kernel$ is \emph{asymptotically smooth}, i.e.,
there exist constants $C,r>0$ such that
\[
  \bigg|\frac{\partial^{|\bs\alpha|+|\bs\beta|}}
  	{\partial{\bs x}^{\bs\alpha}
  	\partial{\bs y}^{\bs\beta}} \kernel ({\bs x},{\bs y})\bigg|
  		\leq C \frac{(|\bs\alpha|+|\bs\beta|)!}
		{r^{|\bs\alpha|+|\bs\beta|}
		\|{\bs x}-{\bs y}\|_2^{|\bs\alpha|+|\bs\beta|}}
\]
for all ${\bs x}\neq
{\bs y}$, ${\bs x},{\bs y}\in\Omega$, {for all} $\bs\alpha,\bs\beta\in\mathbb{N}^d$.
This means that \(\kernel\) is an analytic function whenever
\(\|{\bs x}-{\bs y}\|_2\geq\rho\) for a given \(\rho>0\).
Then, as a consequence of the vanishing moment property,
\({\bs K}^\Sigma\) becomes quasi-sparse and can efficiently be computed.
In detail, there holds the following result, see \cite{HM24}.

\begin{theorem}\label{thm:compression}
Let \(X\) be quasi-uniform with \(\#X=N\) and
set all coefficients of ${\bm K}^\Sigma$ from 
\eqref{eq:K-samplet-basis} to zero which satisfy the admissibility condition
{\begin{equation}\label{eq:admissibility}
     \dist(\tau,\tau')\ge\rho\max\{\diam(\tau),\diam(\tau')\}, \quad \rho>0,
\end{equation}}
where \(\tau\) is the cluster supporting \(\sigma_{j,k}\) and \(\tau'\) is the 
cluster supporting \(\sigma_{j',k'}\), respectively. 
{Then, there exists
a constant \(C>0\), such that the resulting compressed matrix 
${\bs K}^{\Sigma,\rho}$ satisfies}
\begin{equation}\label{eq:consistency_error}
{\big\|{\bs K}^\Sigma-{\bs K}^{\Sigma,\rho}\big\|_F}
   \leq C \bigg(\frac{r\rho}{d}\bigg)^{-2(q+1)}{\big\|{\bs K}^\Sigma\big\|_F}.
\end{equation}
The compressed matrix \({\bs K}^{\Sigma,\rho}\) has $\mathcal{O}(N\log N)$ nonzero
entries.
\end{theorem}

\begin{figure}[htbp!]
    \centering
        \includegraphics[width=0.24\linewidth]{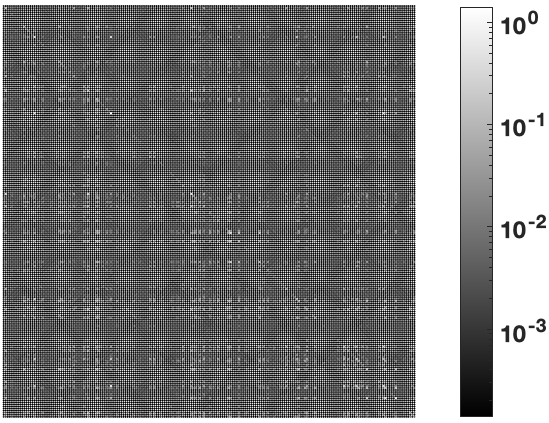} \quad \quad \quad  \quad 
       \includegraphics[width=0.24\linewidth]{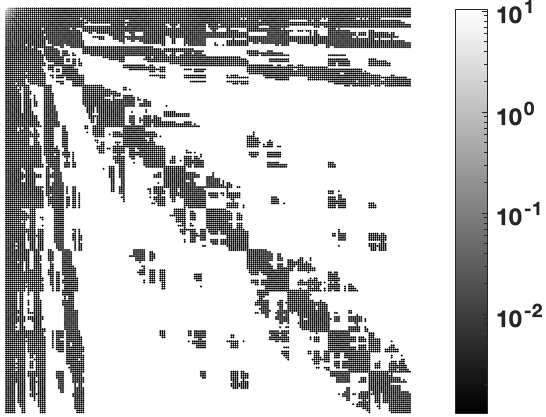}
    \caption{Visualization of the kernel matrices $\bm K^{\Sigma}$ and $\bm K^{\Sigma,\rho}$ 
    (left and right), constructed using the Matérn-$3/2$ kernel on $N = 10\,000$ points 
    uniformly distributed in the unit square $[0,1]^2$. Both matrices are displayed on 
    a logarithmic color scale to reveal the large proportion of small entries, with each 
    matrix subdivided into $200 \times 200$ blocks colored by mean absolute value. 
    The transition from $\bm K^{\Sigma}$ to $\bm K^{\Sigma,\rho}$ is obtained by 
    thresholding small coefficients to zero, reducing the number of nonzero entries 
    to $\mathcal{O}(N \log N)$.}
    \label{fig:pattern_KCompressed}
\end{figure}

Theorem~\ref{thm:compression} establishes that, for asymptotically smooth kernels, 
the kernel matrix admits a compressed, sparse representation in the samplet basis 
with controllable accuracy. 
In contrast, the matrix $\bm K^\Sigma$ exactly represents the kernel matrix
in the samplet basis.
A visualization of the sparsity pattern of the compressed matrix 
$\bm K^{\Sigma,\rho}$ is shown in Figure~\ref{fig:pattern_KCompressed}.
{We remark that the hidden constant in the big-$\mathcal{O}$ notation in Theorem~\ref{thm:compression} involves the
dimension of the polynomial space 
\(\Pcal_q\), which grows exponentially in 
\(d\). The approach is therefore best suited
for moderate dimensions \(d\).}

\subsection{Tree-adaptive sub-sampling}
\label{sec:adaptive_ls}
The formulations \eqref{eq:LSE} and \eqref{eq:system_samplets} employ all
$N$ data sites for the representation of the interpolant. As a consequence,
the solution of the resulting system is no longer computationally feasible in
the regime of very large $N$, without employing additional compression
techniques.
In this section, we develop an adaptive selection strategy to identify
a  subspace of \(\Hcal_X\) that approximates the orthogonal
projection of a given function \(h \in \Hcal\) onto \(\Hcal_X\),
cf.~\eqref{eq:HX}, up to a prescribed consistency error.

To this end, we employ the \emph{adaptive tree approximation}
from \cite{binev2004fast}, where the initial tree is obtained from the cluster
tree \(\Tcal\) underlying the samplet basis.
Then, given a relevant subtree, we select a representative data site from
each remaining cluster to construct the desired subspace.

We start from the samplet transformed data values 
${\bs h}^\Sigma = {\bs T}{\bs h}$, where \({\bs h}=[h({\bs x}_i)]_{i=1}^N\)
and apply tree-coarsening from \cite{binev2004fast} to this setting.
We define the \emph{energy} contained in a cluster \(\tau\in\Tcal\) as the sum of energies of its children and add the 
squared Euclidean norm of the
samplet coefficients belonging to \(\tau\). 
This yields
\begin{equation}\label{eq:energy}
e(\tau)\isdef\sum_{\tau'\in\operatorname{child}(\tau)}e(\tau')
+\sum_{{(j,k) \, : \, \operatorname{supp}\sigma_{j,k}\subset\tau}}
\big({\sigma_{j,k}(h)}\big)^2,
\end{equation}
{where the second sum ranges over the samplets $\sigma_{j,k}$ of the
basis introduced in Section~\ref{sec:samplets} and
$\sigma_{j,k}(h)$ are the corresponding entries of the
transformed data vector ${\bs h}^\Sigma={\bs T}{\bs h}$. For the root cluster,
the sum additionally comprises the coefficients of the scaling distributions
spanning $\Xcal_0'$.} Observe that the quantity \(e(\tau)\) is the contribution
of the subtree rooted in \(\tau\) to the squared Euclidean norm
\(\big\|{\bs h}^\Sigma\big\|_2^2\). Especially, there holds
\(e(X)=\big\|{\bs h}^\Sigma\big\|_2^2\). Due to the orthogonality of the samplet
basis, the energy \(e(X)\) is localized in the clusters \(\tau\in\Tcal\).

Based on the energies \eqref{eq:energy}, we next define
\[
\tilde{e}(\tau') \isdef
\frac{\sum_{\mu\in\operatorname{child}(\tau)}e(\mu)}
{e(\tau)+\tilde{e}(\tau)}\tilde{e}(\tau)
\quad\text{for all}\ \tau'\in\operatorname{child}(\tau) 
\]
with \(\tilde{e}(X)\isdef e(X)\). 
{Note that \(e(\tau)\) is large for
every ancestor of a significant coefficient, whereas \(\tilde e(\tau)\) discounts
the energy already resolved by the ancestors of \(\tau\) and hence reflects its
relevance within the hierarchy. The 
energies \(e(\tau)\) are accumulated
bottom-up, starting from the leaves, whereas the modified energies
\(\tilde e(\tau)\) are propagated top-down, starting from the root.} 
Employing the modified energies
\(\tilde{e}(\tau)\), \(\tau\in\Tcal\), we apply the thresholding version
of the second algorithm from \cite{binev2004fast} with threshold
\(t = \varepsilon^2\|{\bs h}^\Sigma\|_2^2=\varepsilon^2\|{\bs h}\|_2^2\).
The result of this procedure is a subtree \(\Tcal_t\) that approximates
\({\bs h}^\Sigma\) up to a relative error of 
{\(\varepsilon\)} 
in the Euclidean
norm. 

In view of the discussion in the previous subsection, the Euclidean norm
of the coefficient vector \({\bs h}^{\Sigma}\) corresponds to the
\(\Xcal'\)-norm of the functional
\(
h_X'=\sum_{j,k} h_{j,k}^\Sigma\sigma_{j,k}=\sum_{i=1}^N h(\bm x_i)\delta_{\bm x_i}\in\Hcal',
\)
while its \(\Hcal'\)-norm is given by 
{
$
\|h'_X\|_{\mathcal H'}^2=\|Rh'_X\|_{\mathcal H}^2={\bm h}^\intercal{\bm K}{\bm h}=({\bm h}^\Sigma)^\intercal{\bm K}^\Sigma{\bm h}^\Sigma,$
}
where \(R\) denotes the Riesz isometry.
{Equivalently, $\|h'_X\|_{\Hcal'}$ is the native space norm of the
function $Rh_X'=\sum_i h(\bm x_i)\kernel(\bm x_i,\cdot)\in\Hcal_X$. Observe that this is {not} the native space norm of the kernel interpolant of
$h$, which is instead given by ${\bm h}^\intercal{\bm K}^{-1}{\bm h}$.}
Different from the \(\Xcal'\)-norm, the \(\Hcal'\)-norm is
not localized by the samplet coefficients in the nodes of the cluster tree.
In particular, the entries
$k^\Sigma_{(j,k),(j',k')}=\langle\psi_{j,k},\psi_{j',k'}\rangle_{\mathcal{H}}$
of the kernel matrix in samplet coordinates
with \((j,k)\neq(j',k')\) capture cross-scale and spatial interactions.
This also does not change when replacing \({\bs K}^\Sigma\) by its compressed
counterpart \({\bs K}^{\Sigma,\rho}\). 
{
To extend the algorithm of \cite{binev2004fast} to the
$\Hcal'$-setting, we localize an upper bound
of the $\Hcal'$-norm. There holds
\[
({\bs h}^\Sigma)^\intercal{\bs K}^\Sigma{\bs h}^\Sigma
=\sum_{(j,k),(j',k')}
\frac{k^\Sigma_{(j,k),(j',k')}}{\|\psi_{j,k}\|_{\Hcal}\|\psi_{j',k'}\|_{\Hcal}}
\big(\|\psi_{j,k}\|_{\Hcal}h^\Sigma_{j,k}\big)
\big(\|\psi_{j',k'}\|_{\Hcal}h^\Sigma_{j',k'}\big).
\]
Estimating $|ab|\le\frac12(a^2+b^2)$ for the last two factors and using the
symmetry of ${\bs K}^\Sigma$, we obtain
\begin{equation}\label{eq:majorant}
({\bs h}^\Sigma)^\intercal{\bs K}^\Sigma{\bs h}^\Sigma
\le\sum_{(j,k)}\|\psi_{j,k}\|_{\Hcal}^2\big(h^\Sigma_{j,k}\big)^2
\sum_{(j',k')}
\frac{\big|k^\Sigma_{(j,k),(j',k')}\big|}{\|\psi_{j,k}\|_{\Hcal}\|\psi_{j',k'}\|_{\Hcal}}
\defis({\bs h}^\Sigma)^\intercal{\bs D}_w{\bs h}^\Sigma .
\end{equation}
Since $\|\psi_{j,k}\|_{\Hcal}^2=k^\Sigma_{(j,k),(j,k)}$,
the entries of the diagonal matrix
${\bs D}_w$ are given by
\begin{equation}\label{eq:Dw}
d_{j,k}\isdef k^\Sigma_{(j,k),(j,k)}
+\sum_{(j',k')\neq(j,k)}\big|k^\Sigma_{(j,k),(j',k')}\big|
\sqrt{\frac{k^\Sigma_{(j,k),(j,k)}}{k^\Sigma_{(j',k'),(j',k')}}}>0 .
\end{equation}
Replacing the squared coefficients 
in \eqref{eq:energy} by the
ones weighted by ${\bs D}_w$, we define
\begin{equation}\label{eq:energyH}
e_{\Hcal}(\tau)\isdef\sum_{\tau'\in\operatorname{child}(\tau)}e_{\Hcal}(\tau')
+\sum_{(j,k)\,:\,\operatorname{supp}\sigma_{j,k}\subset\tau}
d_{j,k}\big(\sigma_{j,k}(h)\big)^2 .
\end{equation}
Since the energies are nonnegative and additive with
$e_{\Hcal}(X)=({\bs h}^\Sigma)^\intercal{\bs D}_w{\bs h}^\Sigma$,  the
algorithm of \cite{binev2004fast} applies without modifications.
In analogy to the Euclidean case, we set $t=\varepsilon^2e_{\Hcal}(X)$.
Let ${\bs h}^\Sigma_t$ be obtained from ${\bs h}^\Sigma$ by setting to zero the
coefficients of the clusters not contained in $\Tcal_t$, and let $h_{X,t}'$
denote the corresponding functional. Since ${\bs D}_w$ is diagonal,
$({\bs h}^\Sigma-{\bs h}^\Sigma_t)^\intercal{\bs D}_w
({\bs h}^\Sigma-{\bs h}^\Sigma_t)$ is the discarded energy. Hence, if the
discarded energy is at most $t$, applying \eqref{eq:majorant} to
${\bs h}^\Sigma-{\bs h}^\Sigma_t$ yields
\[
\big\|h_X'-h_{X,t}'\big\|_{\Hcal'}^2
=({\bs h}^\Sigma-{\bs h}^\Sigma_t)^\intercal{\bs K}^\Sigma
({\bs h}^\Sigma-{\bs h}^\Sigma_t)
\le({\bs h}^\Sigma-{\bs h}^\Sigma_t)^\intercal{\bs D}_w
({\bs h}^\Sigma-{\bs h}^\Sigma_t)
\le\varepsilon^2({\bs h}^\Sigma)^\intercal{\bs D}_w{\bs h}^\Sigma .
\]
Dividing by $\|h_X'\|_{\Hcal'}^2$, we obtain
\begin{equation}\label{eq:relerrH}
\frac{\|h_X'-h_{X,t}'\|_{\Hcal'}^2}{\|h_X'\|_{\Hcal'}^2}
\le\varepsilon^2\,
\frac{({\bs h}^\Sigma)^\intercal{\bs D}_w{\bs h}^\Sigma}
{({\bs h}^\Sigma)^\intercal{\bs K}^\Sigma{\bs h}^\Sigma}.
\end{equation}
The left-hand side of \eqref{eq:relerrH} is the {squared} relative error in the
$\Hcal'$-norm.
By \eqref{eq:majorant}, the ratio on the right-hand side is at least 1 and
measures by how much the localized energy overestimates the $\Hcal'$-norm. It
serves as an indicator for the quality of the localization. 
If it is close to
1, the energy $e_{\Hcal}(X)$ essentially coincides with
$\|h_X'\|_{\Hcal'}^2$
and $\Tcal_t$ approximates $h_X'$ up to a relative error of essentially
$\varepsilon$ in the $\Hcal'$-norm.}

\begin{figure}[htb]
\centering
    \includegraphics[width=0.2\linewidth]{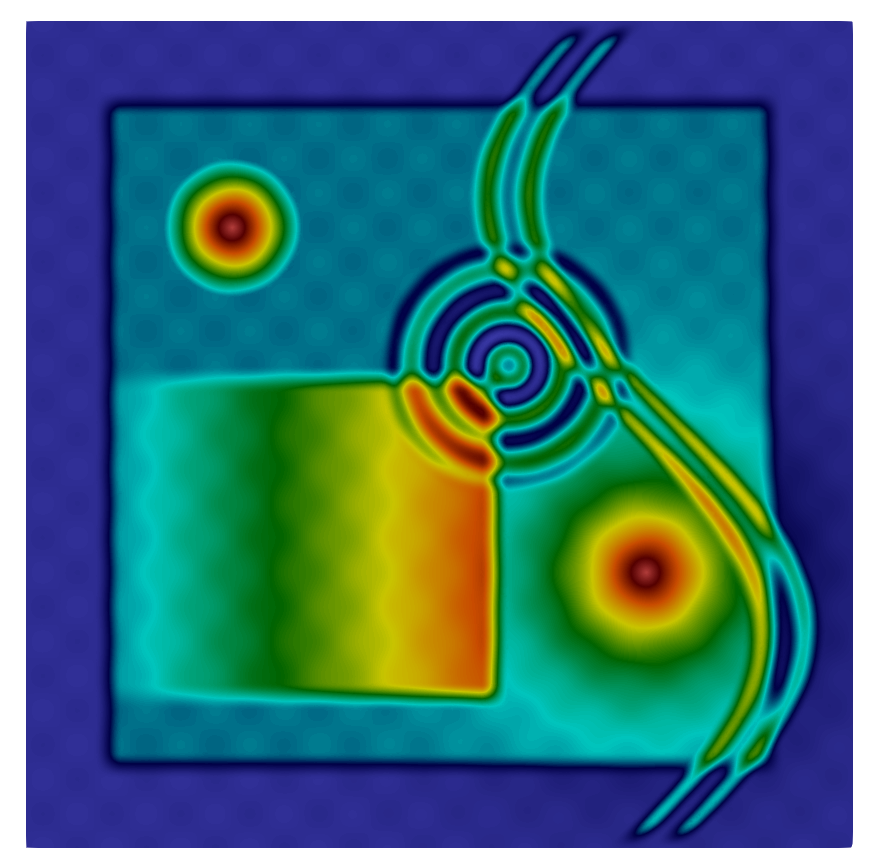} 
    \quad \quad \quad
   \includegraphics[width=0.2\linewidth]{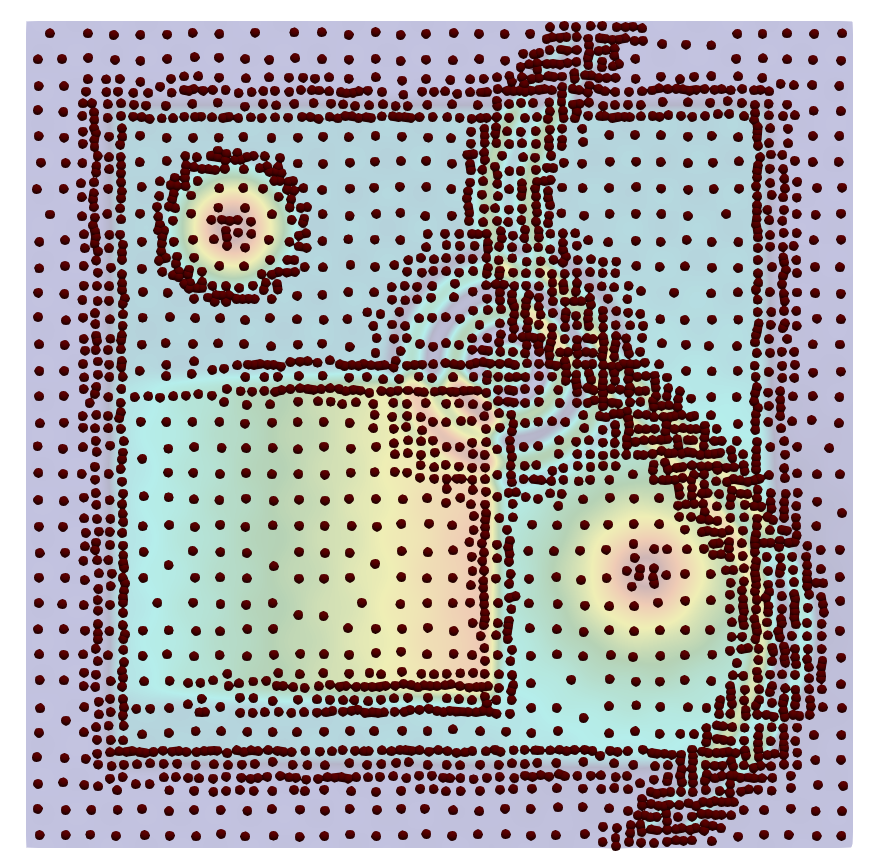} 
   \quad \quad \quad
    \includegraphics[width=0.2\linewidth]{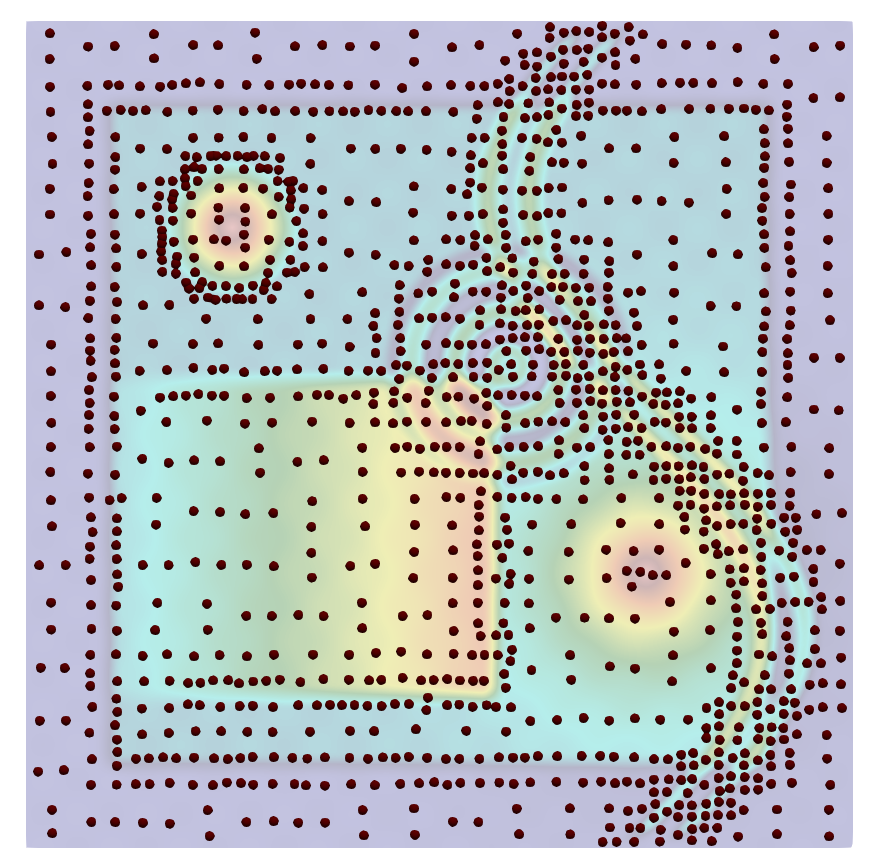} 
\caption{\label{fig:mySubsamples} {Surface plot of the data generating function} (left), adaptive subsample using
the $\Xcal'$-norm (middle) and $\Hcal'$-norm corresponding to a Mat\'ern-1/2
kernel (right). The threshold is set to $\varepsilon^2=10^{-8}$.}
\end{figure}

\begin{remark}
{The Euclidean case is recovered for ${\bs K}^\Sigma={\bs I}$, which
corresponds to the $\ell^2$-kernel, since then ${\bs D}_w={\bs I}$ and the
ratio in \eqref{eq:relerrH} equals 1.} The $\ell^2$-kernel
does not introduce correlations between data sites,
whereas less localized kernels model correlations between them.
Therefore, the choice of the $\Xcal'$-norm for determining a subsample is
the most conservative one, while the choice of the $\Hcal'$-norm of
a less localized kernel yields a sparser subsample. We refer to
Figure~\ref{fig:mySubsamples} for an example.
\end{remark}

In order to quantify the relevance of a cluster, either with respect to the
\(\Xcal'\)-norm or the \(\Hcal'\)-norm, we associate to each $\tau\in\Tcal$
the corresponding energy functional.
To this end, following \cite{binev2004fast}, we proceed in two stages:
\begin{itemize}
  \item In a \emph{bottom-up} traversal of $\Tcal$, we accumulate for each
  cluster $\tau$ the local samplet energy and the contributions from its
  descendants, resulting in a functional $e(\tau)$ that measures the energy of
  the subtree rooted at $\tau$.
  \item In a subsequent \emph{top-down} traversal, we redistribute the global
  energy along the tree to obtain a propagated functional $\tilde e(\tau)$,
  which reflects the relative importance of each cluster within the entire
  hierarchy.
\end{itemize}

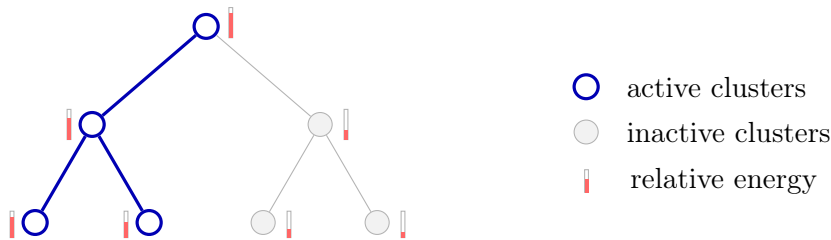
\begin{figure}[htbp]
  \centering
  \begin{tikzpicture}[
    level distance=13mm,
    level 1/.style={sibling distance=30mm},
    level 2/.style={sibling distance=15mm},
    level 3/.style={sibling distance=10mm},
    every node/.style={circle, draw, inner sep=2pt, minimum size=9pt},
    inactive/.style={draw=gray!60, fill=gray!10},
    active/.style={draw=blue!70!black, very thick},
    edge from parent/.style={draw=gray!60},
    active edge/.style={draw=blue!70!black, very thick},
    energyaxis/.style={draw=gray!60, line width=0.4pt},
    energybar/.style={draw=none, fill=red!60}
  ]
  \node[active] (root) {}
    child [edge from parent/.style=active edge] { node[active] (A) {}
      child [edge from parent/.style=active edge] { node[active] (A1) {}
      }
      child [edge from parent/.style=active edge] { node[active] (A2) {}
      }
    }
    child [edge from parent/.style={draw=gray!60}] { node[inactive] (B) {}
      child [edge from parent/.style={draw=gray!60}] { node[inactive] (B1) {}
      }
      child [edge from parent/.style={draw=gray!60}] { node[inactive] (B2) {}
      }
    };
  \draw[energyaxis] ($(root.east) + (0.12,-0.15)$) rectangle ++(0.05,0.40);
  \draw[energybar]  ($(root.east) + (0.12,-0.15)$) rectangle ++(0.05,0.32);
  \draw[energyaxis] ($(A.west) + (-0.15,-0.20)$) rectangle ++(0.05,0.40);
  \draw[energybar]  ($(A.west) + (-0.15,-0.20)$) rectangle ++(0.05,0.29);
  \draw[energyaxis] ($(B.east) + (0.15,-0.20)$) rectangle ++(0.05,0.40);
  \draw[energybar]  ($(B.east) + (0.15,-0.20)$) rectangle ++(0.05,0.13);
  \draw[energyaxis] ($(A1.west) + (-0.15,-0.20)$) rectangle ++(0.05,0.35);
  \draw[energybar]  ($(A1.west) + (-0.15,-0.20)$) rectangle ++(0.05,0.27);
  \draw[energyaxis] ($(A2.west) + (-0.15,-0.20)$) rectangle ++(0.05,0.35);
  \draw[energybar]  ($(A2.west) + (-0.15,-0.20)$) rectangle ++(0.05,0.21);
  \draw[energyaxis] ($(B1.east) + (0.15,-0.20)$) rectangle ++(0.05,0.35);
  \draw[energybar]  ($(B1.east) + (0.15,-0.20)$) rectangle ++(0.05,0.11);
  \draw[energyaxis] ($(B2.east) + (0.15,-0.20)$) rectangle ++(0.05,0.35);
  \draw[energybar]  ($(B2.east) + (0.15,-0.20)$) rectangle ++(0.05,0.07);
  
  \begin{scope}[shift={(3.8,-0.8)}, every node/.style={}]
    \node[active, draw, circle, inner sep=2pt, minimum size=9pt] at (1.2,0) {};
    \node[anchor=west, draw=none] at (1.6,0) {active clusters};
    
    \node[inactive, draw, circle, inner sep=2pt, minimum size=9pt] at (1.2,-0.6) {};
    \node[anchor=west, draw=none] at (1.6,-0.6) {inactive clusters};
    
    \draw[energyaxis] (1.18,-1.4) rectangle ++(0.05,0.30);
    \draw[energybar]  (1.18,-1.4) rectangle ++(0.05,0.18);
    \node[anchor=west, draw=none] at (1.65,-1.25) {relative energy};
  \end{scope}
  
  \end{tikzpicture}
  \caption{Illustration of a cluster tree (gray) and an adaptive subtree (blue)
  selected according to samplet energies. Each node represents a cluster.
  Active clusters belong to the adaptive tree ${\mathcal{T}_t}$
  and the red bar below a node indicates the relative magnitude of its
  associated energy.}
  \label{fig:adaptiveTreeIllustration}
\end{figure}

The clusters are then sorted in decreasing order of $\tilde e(\tau)$ and the
adaptive subtree ${\Tcal}_t\subset\Tcal$ is constructed by adding clusters
until the remaining energy falls below the prescribed threshold.
A simple representative picture illustrating the procedure can be found in
Figure~\ref{fig:adaptiveTreeIllustration}.

In view of Theorem~\ref{theo:waveletProperties}.(iv), 
the adaptive subtree ${\Tcal}_t$ localizes the regions of the data
domain, where the data exhibits singularities,
see also \cite{avesani2025multiresolution}, resulting in a locally deeper tree
with smaller leaves. Thus, to obtain a reduced set of representative data sites,
we associate
to each remaining leaf cluster $\tau\in{\Tcal}_t$ a single
representative index $\iota(\tau)\in\{1,\dots,N\}$. In practice, $\iota(\tau)$
is chosen as the index of the data site in $\tau$ that is closest to the
cluster's centroid. Another possible choice is to randomly select centers by a
suitably defined density, see \cite{HM24}.
Associated to the set of representative indices
\(
I
\isdef
\{ \iota(\tau) : \tau \text{ is a leaf of } {\Tcal}_t\}
\subset \{1,\dots,N\},
\) 
we introduce the set representative data sites
\(
{X}_t
\isdef
\{\bm x_i : i \in I\},
\) which reflects the energy distribution of the data.
Clusters carrying negligible energy are discarded, whereas clusters with
significant contributions are represented by the cluster representative center.


\section{Multiscale kernel lasso}
\label{sec:newton_method}
In this section, we introduce an $\ell^1$-regularized formulation of the
kernel approximation problem, which promotes sparse representations in
the multiscale setting, see, e.g., \cite{Candes,Donoho,daubechies2004iterative}.
We first describe a kernel basis pursuit formulation,
and then present the semismooth Newton method used to solve the resulting
optimization problem, including its efficient realization based on an
online SVD and a trust-region strategy.

\subsection{Kernel basis pursuit}\label{sec:kernel_basis_pursuit}
{Beyond providing a computationally efficient solver for large-scale
kernel approximation problems, the formulation developed in this section
also addresses a second, complementary goal: selecting an appropriate
kernel (or an appropriate combination of kernel type and lengthscale
parameter) for a given data set, in situations where little a priori
information on a suitable kernel is available. Rather than fixing a single
kernel $\kernel$ beforehand, we consider a \emph{dictionary} of candidate
kernels $\kernel_1,\ldots,\kernel_L$, differing, e.g., in kernel type or in
lengthscale, and let a subsequent $\ell^1$-regularization automatically
select the relevant kernel contributions from the dictionary, effectively
performing kernel and lengthscale selection jointly with the
approximation.}

We consider the sets of data sites
\(X_\ell\isdef\big\{{\bs x}_1^{(\ell)},\ldots,{\bs x}_{N_\ell}^{(\ell)}\big\}\),
\(\ell=1,\ldots,L\), and
a corresponding dictionary of kernels 
 $\kernel_1,\ldots,\kernel_L$ and seek a (sparse) representation of the form
	\begin{equation}\label{eq:mskernel}
	{s_h^{\rm ms}}=\sum_{\ell=1}^L\sum_{i=1}^{N_\ell} \alpha_{\ell,i}
	\kernel_\ell\big(\cdot,{\bs x}_i^{(\ell)}\big).
	\end{equation}
{Here, each $X_\ell$ is the set of centers associated to the $\ell$-th
kernel in the dictionary, and different kernels may in principle be
evaluated at different, possibly overlapping, sets of data sites. In
particular, in the setting of interest here, each
$X_\ell$ can be chosen as a copy of the reduced, tree-adaptively selected
data site set $X_t$ from Section~\ref{sec:adaptive_ls}, equipped with the
$\ell$-th kernel of the dictionary; this is precisely the setting used in
our numerical experiments in Section~\ref{sec:results}.}
This representation particularly comprises the multiscale scattered data
interpolation framework, see for example \cite{CDS98,MZ,tropp2004greed},
if a sequence of kernels with decreasing lengthscale
parameters is chosen associated to a sequence of nested point sets.
Setting \(M\isdef N_1+\ldots+N_L\) as well as
	\[
	{\bs K}\isdef[{\bs K}_1,\ldots,{\bs K}_L]\in\Rbb^{N\times M},\quad 
	{\bs K}_\ell\isdef\big[
	\kernel_\ell\big({\bs x}_i,{\bs x}_j^{(\ell)}\big)
	\big]_{\substack{i=1,\ldots,N \\ j=1,\ldots,N_\ell}},
	\]
	and
	\[
	{\bs\alpha}^\intercal
	\isdef[{\bs\alpha}_{1}^\intercal,\ldots,{\bs\alpha}_{L}^\intercal]\in\Rbb^M,
  \quad {\bs\alpha}_\ell\in\Rbb^{N_\ell},
	\]
it is well known that a sparse solution to \eqref{eq:mskernel} is obtained by
solving the convex optimization problem
\begin{equation}\label{eq:SSLASSO}
	\min_{{\bs\alpha}\in\Rbb^{M}}\;
	\frac{1}{2}\|{\bs h}-{\bs K}{\bs\alpha}\|^2_2
	+\sum_{i=1}^{M} w_i|\alpha_i|,
\end{equation}
for a given weight vector ${\bs w}\in\Rbb^{M}$. The entries of 
${\bs w}=[w_i]_{i=1}^M$ act as
component-wise regularization parameters, thereby controlling the level of
sparsity in each coefficient individually. 
We refer to, e.g., \cite{daubechies2004iterative,Lorenz,RT} 
for a detailed analysis of the regularizing properties
and strategies for appropriate parameter selection.

In contrast to the single scale formulation \eqref{eq:SSLASSO}, it has been
observed in \cite{BaroliHarbrechtMulterer2024} that a richer class of sparse
representations can be obtained by employing a samplet basis. This yields the
optimization problem
	\begin{equation}\label{eq:MSLASSO}
		\min_{{\bs\beta}\in\Rbb^M}
		\frac{1}{2}\|{\bm T \bs h}-{\bs K^\Sigma}{\bs\beta}\|^2_2
		+\sum_{i=1}^M w_i|\beta_i|
        \end{equation}
with
\[
	{\bs K}^\Sigma\isdef[{\bs T}{\bs K}_1{\bs T}_1^\intercal
	,\ldots,{\bs T}{\bs K}_L{\bs T}_{L}^\intercal]\in\Rbb^{N\times M}
	\]
	and
	\[
	{\bs\beta}^\intercal
	\isdef[{\bs\beta}_{1}^\intercal,\ldots,{\bs\beta}_{L}^\intercal]\in\Rbb^M,
  \quad 
	{\bs\beta}_\ell\isdef{\bs T}_\ell{\bs\alpha}_\ell\in\Rbb^{N_\ell},
	\]
where \({\bs T}_\ell\) denotes the samplet transform on
\(X_\ell\).
Due to the orthonormality of the samplet basis, there particularly holds
\(
\|{\bs h}-{\bs K}{\bs\alpha}\|_2=
\|{\bm T \bs h}-{\bs K^\Sigma}{\bs\beta}\|_2.
\)
Hence, the formulations \eqref{eq:SSLASSO} and \eqref{eq:MSLASSO} only differ
in the choice of the regularization term.

\subsection{Semismooth Newton method}\label{sec:SSN}
To numerically solve \eqref{eq:SSLASSO} and \eqref{eq:MSLASSO}, we apply the
trust-region semismooth Newton approach proposed in \cite{ouyang2025trust}. 
To compute the Newton update in a numerically stable way, we employ the online
SVD approach from \cite{brand2002incremental}. For the reader's convenience, we
recall the key algorithms and adapt them to our framework. To this end, 
we start from the minimization problem
\begin{align}\label{eq:Composite}
	\min_{{\bs\alpha}\in\Rbb^M}\; \Phi({\bs\alpha}),\quad\text{where } 
	\Phi({\bs\alpha})\isdef f({\bs\alpha}) + g({\bs\alpha})
\end{align}
with
\[
f({\bs\alpha}) \isdef \tfrac{1}{2}\|{\bs h}-{\bs K}{\bs\alpha}\|_2^2,
\quad
g({\bs\alpha}) \isdef \sum_{i=1}^M w_i|\alpha_i|.
\]
The functional $f\colon\Rbb^M\to\Rbb$ is continuously differentiable
and reflects the data-fidelity, while $g\colon\Rbb^M\to[0,\infty)$
is a convex, lower semicontinuous and proper mapping that enforces sparsity.
The vector \({\bs\alpha}\in\Rbb^M\) is a critical point of $\Phi$ if and only
if
\[
\Fnat{{\bs\alpha}}\isdef
{\bs\alpha} - \operatorname{SS}_{\lambda \bs w}\big({{\bs\alpha} 
- \lambda \nabla f({\bs\alpha})\big)} = \bs 0.
\]
Herein, the \emph{soft-shrinkage operator} is given by
\[
[\operatorname{SS}_{\lambda \bs w}(\bs \alpha)]_i
\isdef \operatorname{sign}(\alpha_i)\,\max(|\alpha_i|-\lambda w_i,0),
\quad i=1,\ldots,M.
\]

Following the normal-map reformulation from \cite{ouyang2025trust},
we introduce the associated normal-map operator
\begin{equation*}
\begin{aligned}
	\Fnor{{\bs z}} 
	&\isdef \nabla f\big(\operatorname{SS}_{\lambda \bs w}({{\bs z}})\big)
	+\frac{1}{\lambda}\big({\bs z} 
  - \operatorname{SS}_{\lambda \bs w}({{\bs z}})\big)\\
&={\bs K}^\intercal\big({\bs K}\operatorname{SS}_{\lambda{\bs w}}({\bs z}) 
  - {\bs h}\big) +\frac{1}{\lambda}\big({\bs z}
-\operatorname{SS}_{\lambda{\bs w}}({\bs z})\big),
\end{aligned}
\end{equation*}
If ${\bs z}\in\Rbb^M$ satisfies $\Fnor{{\bs z}}=\mathbf 0$, then ${\bs\alpha}
=\operatorname{SS}_{\lambda \bs w}({{\bs z}})$ is a critical point of $\Phi$. 
Conversely, if ${\bs\alpha}$ is a critical point of $\Phi$,
then ${\bs z}={\bs\alpha}-\lambda\nabla f({\bs\alpha})$ satisfies 
$\Fnor{{\bs z}}=\mathbf 0$.

To compute the roots of the normal map, we employ the semismooth Newton method
\(
	\bs z_{k+1} = \bs z_k - {\bs M}_k^{-1} \Fnor{\bs z_k},
\)
where
$
{\bs M}_k\isdef{\bs K}^\intercal {\bs K}{\bs I}_{\Acal_k} 
+ \frac{1}{\lambda}\mathbf{I}_{\Ical_k}.
$
Herein, we define the active and inactive sets as
	\begin{equation}\label{eq:act_inact}
		\Act_k \isdef \{i \in \{1,\dots,M\} : |z_{k,i}| > \lambda w_i\}, \quad 
		\Inact_k \isdef \{i \in \{1,\dots,M\}:|z_{k,i}| \leq \lambda w_i\}
	\end{equation}
and \({\bs I}_{\Acal_k}\in\Rbb^{M\times M}\),
\({\bs I}_{\Ical_k}\in\Rbb^{M\times M}\)
are partial identity matrices corresponding to the active- and inactive indices,
respectively.
{Observe that
${\bs I}_{\Acal_k}$ and ${\bs I}_{\Inact_k}={\bs I}-{\bs I}_{\Acal_k}$ are
complementary orthogonal projectors, i.e.\
${\bs I}_{\Acal_k}{\bs I}_{\Inact_k}=\bs 0$ and
${\bs I}_{\Acal_k}^\intercal={\bs I}_{\Acal_k}$.
The full Newton step $\bs s_k\isdef -{\bs M}_k^{-1}F_{\mathrm{nor}}^\lambda(\bs z_k)$
solves
\[
{\bs M}_k\bs s_k = -F_{\mathrm{nor}}^\lambda(\bs z_k).
\]
Multiplying this equation from the left by ${\bs I}_{\Acal_k}$ and using
${\bs I}_{\Acal_k}{\bs I}_{\Inact_k}=\bs 0$ gives
\[
{\bs I}_{\Acal_k}{\bs M}_k
={\bs I}_{\Acal_k}\Big({\bs K}^\intercal {\bs K}{\bs I}_{\Acal_k} 
+ \frac{1}{\lambda}\mathbf{I}_{\Ical_k}\Big)
={\bs I}_{\Acal_k}{\bs K}^\intercal {\bs K}{\bs I}_{\Acal_k},
\]
since the second term is annihilated by the projection.
Hence the active block of the Newton step satisfies the symmetric,
positive semidefinite reduced system
\begin{equation}\label{eq:ssn_active}
	{\bs I}_{\Acal_k}{\bs K}^\intercal {\bs K}{\bs I}_{\Acal_k} \, \bs s_k
		= - {\bs I}_{\Acal_k}F_{\mathrm{nor}}^\lambda({\bs z}_k).
\end{equation}
Note that both sides of
\eqref{eq:ssn_active} vanish identically on the inactive indices
$\Inact_k$, so \eqref{eq:ssn_active} only determines the components
$[\bs s_k]_i$, $i\in\Acal_k$; the complementary inactive components of the
lifted step are instead recovered explicitly from the trust-region approach
given further below.}

\subsubsection*{Online SVD}
A numerically stable solver for \eqref{eq:ssn_active} is obtained from
the online SVD of \({\bs K}{\bs I}_{\Acal_k}\), which keeps track of the columns
of \({\bs K}\) that have been active up to iteration \(k\). A block version of
the corresponding procedure from \cite{brand2002incremental} is shown in
Algorithm~\ref{alg:onlineSVD_matrix}.

\begin{algorithm}[htb]
\caption{Block online SVD}
\label{alg:onlineSVD_matrix}
  \begin{algorithmic}[1]
    \Require Current SVD factors 
      $\bm U \in \mathbb{R}^{N \times \ell}$, 
      $\bm S \in \mathbb{R}^{\ell \times \ell}$, 
      $\bm V \in \mathbb{R}^{s \times \ell}$,
      new columns $\bm K \in \mathbb{R}^{N \times p}$
    \Ensure Updated SVD factors 
      $\bm U\in \mathbb{R}^{N \times (\ell+p)},\ 
      \bm S\in \mathbb{R}^{(\ell+p) \times (\ell+p)},\
       \bm V\in \mathbb{R}^{(s+p) \times (\ell+p)}$
    \State set $\bm C \isdef \bm U^\intercal \bm K$
    \State compute the QR decomposition $ (\bm K - \bm U \bm C) = \bm Q \bm R$
    \State set $\bm U \isdef \begin{bmatrix} \bm U & 
    \bm Q \end{bmatrix}$, $\bm S \isdef 
      \begin{bmatrix}
        \bm S & \bm C \\
        \bm 0 & \bm R
      \end{bmatrix}$, $\bm V \isdef
      \begin{bmatrix}
        \bm V & \bm 0 \\
        \bm 0 & \bm I
      \end{bmatrix}$
    \State compute the SVD 
      $\bm S = \tilde{\bm U} \tilde{\bm S} \tilde{\bm V}^\intercal$
    \State set $\bm U \isdef \bm U \tilde{\bm U}$, $\bm S \isdef \tilde{\bm S}$, 
    $\bm V \isdef \bm V \tilde{\bm V}$
  \end{algorithmic}
\end{algorithm}

The cost of the block online SVD is clearly \(\Ocal\big(N(\ell+p)^2\big)\) if 
\(s\) existing columns are updated by \(p\) new ones and the original rank is
\(\ell\). 
Hence, this approach is particularly useful, if the size of the active set
stays moderate, i.e., \(\ell+p\leq s+p\ll N\).
The stability can be improved
by performing a reorthogonalization step in the computation of the 
expansion coefficients with respect to the current basis of the image in
line 1. In contrast, the orthogonalization of the added vectors by a Householder
QR decomposition is always backward stable.

Now let \(\Dcal\isdef\bigcup_{i=1}^k\Acal_i\) and consider the SVD
\({\bs K}{\bs I}_{\Dcal}={\bs U}{\bs S}{\bs V}^\intercal\). Then, there holds
\(
{\bs I}_{\Dcal}{\bs K}^\intercal{\bs K}{\bs I}_{\Dcal}
={\bs V}{\bs S}^2{\bs V}^\intercal.
\)
From this, we can derive a direct solver for the current active set by computing
the SVD of \({\bs I}_{\Acal_k}{\bs V}{\bs S}\) with cost
\(\Ocal\big((\max\{\ell,\#\Acal_k\})(\min\{\ell,\#\Acal_k\})^2\big)\)
if the rank of the stored online SVD is \(\ell\).

\subsubsection*{Trust-region step management}
{Before detailing the specific construction, let us recall the general
idea underlying trust-region methods.
{At each iteration $k$, a local quadratic model $m_k$ of the objective
is constructed around the current iterate and is trusted only within a ball of radius
$\Delta_k$, the \emph{trust region}. A trial step $\bs s_k$ is computed by
approximately minimizing $m_k$ subject to $\|\bs s_k\|\le \Delta_k$, and
its quality is assessed a posteriori by comparing the actual reduction it
produces with the reduction predicted by the model: if the two agree well,
the step is accepted and the trust region is enlarged; otherwise, the step
is rejected and the trust region is shrunk.}
This mechanism provides a
robust safeguard for the semismooth Newton iteration, in particular in the
presence of active-set changes, where the local linearization underlying
the Newton step may be inaccurate.}
To globalize the semismooth Newton iteration and ensure robustness in the
presence of active-set changes, we employ the trust-region strategy proposed
in \cite{ouyang2025trust} and adapt it to our framework.
At iteration $k$, the reduced trust-region model is constructed on the active
index set $\Act_k$. The corresponding quadratic model reads
\begin{equation}\label{eq:TR_subprob_red}
m_k(\bs q)
=
\iprod{{\bs I}_{\Acal_k}F_{\mathrm{nor}}^\lambda({\bs z}_k)}{\bs q}
+
\frac12
\iprod{{\bs I}_{\Acal_k}{\bs K}^\intercal {\bs K}{\bs I}_{\Acal_k}\bs q}{\bs q}.
\end{equation}
{This model is precisely the restriction, to the active subspace, of a
second-order Taylor-type expansion of the auxiliary function
$\bs z\mapsto\Phi(\operatorname{SS}_{\lambda\bs w}(\bs z))$ around $\bs z_k$.}
Let $\bs q_k$ denote an approximate minimizer of \eqref{eq:TR_subprob_red}
within the trust-region of radius $\Delta_k$. We compute $\bs q_k$ by the
previously described online SVD solver.
The candidate Newton step is subsequently obtained through the lifting procedure
\[
\bar{\bs s}_k
=
\bs q_k
-
\lambda\big(F_{\mathrm{nor}}^\lambda(\bs z_k)+{\bs M}_k \bs q_k\big),
\qquad
\bs s_k
=
\min\bigg\{1,\frac{\Delta_k}{\|\bar{\bs s}_k\|_2}\bigg\}\bar{\bs s}_k,
\]
which yields an approximate semismooth Newton direction satisfying the
trust-region constraint, see \cite{ouyang2025trust}.
{Intuitively, the lifting step translates the (reduced, active-set)
model minimizer $\bs q_k$ back into a genuine approximate solution of the
full semismooth Newton system for $F_{\mathrm{nor}}^\lambda$, and the
subsequent rescaling by $\min\{1,\Delta_k/\|\bar{\bs s}_k\|_2\}$ guarantees
that the resulting step $\bs s_k$ itself respects the trust-region radius
$\Delta_k$, which is not automatically implied by $\|\bs q_k\|\le\Delta_k$
alone once the lifting is applied.}
Due to the active-inactive decomposition, the lifted step admits a simple
structure. For indices $i\in\Act_k$ one has $\bar s_{k,i}=q_{k,i}$, whereas
for $i\in\Inact_k$ the update is given by
\(
\bar s_{k,i}
=
-
\lambda
\big(
\big[F_{\mathrm{nor}}^\lambda(\bs z_k)\big]_i
+
\left[{\bs I}_{\Inact_k}{\bs K}^\intercal
{\bs K}{\bs I}_{\Acal_k}\bs q_k\right]_i\big).
\)
Hence, the Newton update is effectively determined by the active subspace.

Now, step acceptance is based on the standard ratio between actual and predicted
reduction,
$
\tilde{\rho}_k=\frac{\mathrm{ared}_k}{\mathrm{pred}_k}.
$
{Given two thresholds
$0<\eta_1\le\eta_2<1$, the step is classified as \emph{unsuccessful} if
$\tilde\rho_k<\eta_1$: in this case the iterate is not updated,
$\bs z_{k+1}=\bs z_k$, and the trust-region radius is contracted, so that
the next trial step is computed on a smaller, more conservative region. If
instead $\tilde\rho_k\ge\eta_1$, the step is accepted,
$\bs z_{k+1}=\bs z_k+\bs s_k$; in particular, if $\tilde\rho_k\ge\eta_2$ the
agreement between actual and predicted reduction is considered very good
and the trust-region radius is enlarged to allow for longer steps in the
following iteration, while for $\eta_1\le\tilde\rho_k<\eta_2$ the step is
still accepted but the radius is left unchanged.}
Following \cite{ouyang2025trust}, the actual reduction is computed through the
merit function
\begin{equation}\label{eq:merit_fct}
\mer(\bs z)
=
\Phi(\operatorname{SS}_{\lambda \bs w}(\bs z))
+
\frac{\theta\lambda}{2}\|F_{\mathrm{nor}}^\lambda(\bs z)\|_2^2, 
\quad \mathrm{ared}_k=\mer(\bs z_k)-\mer(\bs z_k+\bs s_k),
\end{equation}
with $\theta \in (0,1)$. {We emphasize that a naive, more classical
choice for the acceptance test would instead directly compare the objective
decrease $\Phi(\operatorname{SS}_{\lambda\bs w}(\bs z_k))-\Phi(\operatorname{SS}_{\lambda\bs w}(\bs z_k+\bs s_k))$
with the model decrease $-m_k(\bs q_k)$. However, such a ratio only measures the quality of
the reduced-model step $\bs q_k$ and not of the full lifted step $\bs s_k$
actually taken. The merit function $H_\theta$ instead combines the
value of $\Phi\circ\operatorname{SS}_{\lambda\bs w}$ with the normal-map
residual itself, which allows one to measure the actual reduction achieved
by $\bs s_k$ directly and, together with the matching predicted-reduction
term $\mathrm{pred}_k$ below, is what makes it possible to guarantee that
accumulation points of the resulting iterates are zeros of
$F_{\mathrm{nor}}^\lambda$ and to obtain the fast local convergence
properties established in \cite{ouyang2025trust}.}
The predicted reduction is defined as
\begin{equation}\label{eq:pred_final}
\begin{aligned}
\mathrm{pred}_k
&=
\frac{\theta\|F_{\mathrm{nor}}^\lambda(\bs z_k)\|_2}{2}
\min\{\lambda,\Delta_k,\lambda\|F_{\mathrm{nor}}^\lambda(\bs z_k)\|_2\}
\\
&\quad+
\frac{\nu_k\|F_{\mathrm{nor}}^\lambda(\bs z_k)\|_2}
{\min\{\Delta_k,\lambda\|F_{\mathrm{nor}}^\lambda(\bs z_k)\|_2\}}
\|\operatorname{SS}_{\lambda\bs w}(\bs z_k+\bs s_k)
-
\operatorname{SS}_{\lambda\bs w}(\bs z_k)\|_2^2,
\end{aligned}
\end{equation}
with $\nu_k\in[0,1)$.

The trust-region radius $\Delta_k$ is updated in the standard way based
on the reduction ratio $\tilde{\rho}_k$, also enforcing a lower and upper bound
for the updated radius.

Finally, since the trust-region iteration is formulated in terms of the
normal-map variable $\bs z_k$, we set
$
\bs\alpha_k=\operatorname{SS}_{\lambda\bs w}(\bs z_k).
$
The trust-region semismooth Newton procedure is summarized 
in Algorithm~\ref{alg:trustSSN}.

\begin{algorithm}[htb]
\caption{TR-SSN (Trust-region normal-map semismooth Newton)}
\label{alg:trustSSN}
\begin{algorithmic}[1]
  \Require $\bs K\in \mathbb{R}^{N\times M}$, $\bs h\in \mathbb{R}^N$, 
           $\bs z_0 \in \mathbb{R}^M$, $\bs w \in \mathbb{R}^M$, 
           parameters $\lambda >0$, $0<\eta_1\le \eta_2<1$,
            $0 < \Delta_{\min}\leq \Delta_0 \le \Delta_{\max}$, 
           tolerance $\mathtt{tol} > 0$, maximum iterations
           $\mathtt{maxit} \in \mathbb{N}$
  \Ensure approximate solution \({\bs\alpha}\in\Rbb^M\) to \eqref{eq:Composite}
  \State set $k \isdef 0$, $\Delta_k \isdef \Delta_0$
  \While{$\|F_{\mathrm{nor}}^\lambda(\bs z_k)\|_2 > \mathtt{tol}$ 
  \textbf{and} $k < \mathtt{maxit}$}
    \State compute active and inactive sets $\Act_k,\Inact_k$ as 
    in~\eqref{eq:act_inact}
    \State solve the reduced Newton system for $\bs q_k$ using 
    Algorithm~\ref{alg:onlineSVD_matrix}
 \State set the lifted step 
\[
\bar{\bs s}_k\isdef\begin{cases}
  q_{k,i},&\text{if }i\in\Act_k,\\
  -
\lambda
\left(
\left[F_{\mathrm{nor}}^\lambda(\bs z_k)\right]_i
+
\left[{\bs I}_{\Inact_k}{\bs K}^\intercal
{\bs K}{\bs I}_{\Acal_k}\bs q_k\right]_i
\right), & \text{if }i\in\Inact_k
\end{cases}
\]
\State rescale the step
\[
\bs s_k
\isdef
\min\bigg\{1,\frac{\Delta_k}{\|\bar{\bs s}_k\|_2}\bigg\}\bar{\bs s}_k
\]
\State compute the reduction ratio $\tilde{\rho}_k
\isdef {\mathrm{ared}_k}/{\mathrm{pred}_k}$
according to~\eqref{eq:merit_fct}-\eqref{eq:pred_final}
    \If{$\tilde{\rho}_k < \eta_1$}
      \State \textit{unsuccessful step}
      \State $\bs z_{k+1} \isdef \bs z_k$
      \State $\Delta_{k+1} \isdef \max\{\Delta_{\min},\, \Delta_k/2\}$
    \Else
      \State \textit{successful step}
      \State $\bs z_{k+1} \isdef \bs z_k + \bs s_k$
      \If{$\tilde{\rho}_k \ge \eta_2$}
        \State $\Delta_{k+1} \isdef \min\{\Delta_{\max},\, 2\Delta_k\}$
      \Else
        \State $\Delta_{k+1} \isdef \Delta_k$
      \EndIf
    \EndIf
    \State $k \isdef k+1$
  \EndWhile
  \State \Return $\bs \alpha \isdef \operatorname{SS}_{\lambda{\bs w}}(\bs z_k)$
\end{algorithmic}
\end{algorithm}

\section{Numerical results}\label{sec:results}

In this section, we assess the performance of the proposed framework for
large-scale kernel approximation. In Section~\ref{sec:adaptive_sub_validation},
we first isolate and validate the tree-adaptive subsampling strategy. 
In Section~\ref{sec:TR_validation}, we consider a second set of experiments
with the full pipeline, combining tree-adaptive subsampling, samplet compression, and $\ell^1$-regularized regression solved by the TR-SSN method.
All computations have been performed at the Centro Svizzero di Calcolo
Scientifico (CSCS) on a single node of the Alps cluster\footnote{%
\url{https://www.cscs.ch/computers/alps}}
 with two AMD EPYC 7742
@2.25GHz CPUs with 470GB of main memory. For the computations, we have
used up to 12 cores. The implementation is publicly available at
\url{https://github.com/muchip/fmca}.
All experiments share the same algorithmic framework and compression strategy,
as well as certain parameter choices that are common across tests,
which we describe next.

\paragraph{Samplet compression}
Throughout all experiments, we use the samplet-compressed representation of the
kernel matrix as outlined in Section~\ref{sec:samplets}. 
The parameters for the samplet compression are $\rho$ and $q$.
The first one is associated with the admissibility condition
\eqref{eq:admissibility}, which determines the matrix blocks that are initially truncated.
A larger $\rho$ results in a larger distance between clusters before they
are admissible for compression. In our analysis, we adopt $\rho = d$, so that 
the bound in Equation \eqref{eq:consistency_error} is exponentially governed by 
the number of vanishing moments, which we set $q+1 = 4$.
Furthermore, we apply an a-posteriori thresholding of small matrix entries with
a relative threshold of $\kappa=10^{-7}$ with respect to the Frobenius norm of the matrix.

\paragraph{Choice of the $\ell^1$-weights}
The weights $\bs w$ in the $\ell^1$ penalty are scaled with respect to the
number of selected data sites in order to ensure a consistent balance between
the data-fidelity term and the sparsity regularization.
Let $M$ denote the total number of columns of the kernel matrix $\bm K$.
Then, letting ${X}_t$ be the set of centers obtained by adaptive subsampling,
we set
$
w_i = 1/\sqrt{|{X}_t|}, \ i = 1,\dots,M.
$

\paragraph{Choice of the parameter $\lambda$}
The parameter $\lambda$ defines the normal map $F_{\mathrm{nor}}^\lambda$ and
plays a role analogous to a stepsize in proximal-gradient methods.
Several works based on the same fixed-point mapping suggest choosing $\lambda$
sufficiently small and inversely proportional to the Lipschitz constant of
$\nabla f$, see e.g.~\cite{Xiao2018RegularizedSSN,Davis2017FasterPRADMM,
Zhao2001Monotonicity}. In particular,
if $f$ is twice continuously differentiable, the fixed-point residual is
monotone provided that
$
    0 < \lambda \le \frac{2}{\mathcal{L}},
$
where $\mathcal{L}$ denotes the Lipschitz constant of $\nabla f$.
Following this guideline, we set
$
\lambda = \frac{2}{\mathcal{L}},
$
where $\mathcal{L}$ is approximated in practice via a power iteration applied to
the compressed kernel matrix.

\paragraph{Trust-region radius update}
We adopt the {TR-SSN} algorithm described in
Algorithm~\ref{alg:trustSSN}.
Within this scheme, the trust-region radius $\Delta_k$ is initialized
with $\Delta_0=1$ and updated based on the ratio $\tilde{\rho}_k$
between actual and predicted reduction.
We employ acceptance thresholds $\eta_1=10^{-3}$ and $\eta_2=10^{-1}$
and enforce the bounds
$\Delta_{\min}=10^{-5}$ and $\Delta_{\max}=10^{5}$.

\paragraph{Parameters in the reduction ratio $\tilde{\rho}_k$}
The reduction ratio $\tilde{\rho}_k$ defined in
\eqref{eq:merit_fct}-\eqref{eq:pred_final}
depends on the parameters $\theta$ and $\nu_k$.
These parameters are chosen in accordance with the global and local
convergence conditions established in \cite{ouyang2025trust}
for trust-region semismooth Newton methods applied to the normal-map formulation.
In particular, global convergence requires the parameters $\theta$ and $\nu_k$
to satisfy conditions
involving the Lipschitz constant $\mathcal{L}$ of $\nabla f$.
We therefore set
\[
\theta = \frac{2c_\theta}{\mathcal{L}^2\lambda^2+2} = \frac{c_\theta}{3},
\qquad
\nu = \frac12\min\!\left\{
\theta,\;
c_{\nu}
\left[
1-\frac{\theta}{2}\left(\frac{\mathcal{L}^2\lambda^2}{2}+1\right)
\right]
\right\} = \frac{\theta}{2}
,
\]
with $c_\theta = c_\nu = 0.05$.
The safeguard parameter $\nu_k$ appearing in the predicted reduction
is updated adaptively in order to control the inexactness of the Newton steps.
We define
\begin{equation}\label{eq:nu_k}
\nu_k\isdef \min\Bigl\{\nu, {10^{-3}}\Bigl(
n_{\mathcal S}(k)\log^{2}\big(n_{\mathcal S}(k)\big)
\|\operatorname{SS}_{\lambda\bs w}(\bs z_k+\bs s_k)
-
\operatorname{SS}_{\lambda\bs w}(\bs z_k)\|_2 \Bigr)^{2\tilde p}\Bigr\},
\end{equation}
with $\tilde p = 0.1$ and where $n_{\mathcal S}(k)$ denotes the number of
successful trust-region steps up to iteration $k$.
With this choice, the global convergence conditions remain valid,
while the inexactness tolerance is progressively tightened,
which yields fast local convergence of the TR-SSN method.
In particular, the convergence results of
\cite[Theorem~4.8 and Theorem~6.2]{ouyang2025trust}
ensure global convergence of the normal-map residual and
eventually superlinear (quadratic in our quadratic-$\ell^1$-setting)
local convergence.

\paragraph{Continuation strategy}
To enhance robustness and support sparsity identification,
we adopt a continuation strategy in which the TR-SSN solver 
is applied within an outer loop that progressively reduces the 
$\ell^1$-weight.
Starting from an initial ramp parameter $r_0$, and a given 
$\gamma \in (0,1)$, we define the sequence
\begin{equation}\label{eq:regularized_alg}
    r_j = r_0 \, \gamma^{j}, \qquad j = 0,1,2,\dots, \quad
\end{equation}
and the process is terminated when
$
    r_j < r_{\min},
$
for a prescribed threshold $r_{\min} > 0$.
At each continuation level $j$, we solve the weighted $\ell^1$-regularized 
least-squares problem
\begin{equation}\label{eq:l1_ramp}
\min_{\bm\alpha\in\Rbb^M}\;
\frac12\|\bm h-\bm K \bm\alpha\|_2^2
+ r_j\sum_{i=1}^M w_i |\alpha_i|,
\end{equation}
using TR-SSN with the solution $\bm\alpha_{j-1}$ from the previous level as 
initial guess.
Depending on the specific test case, this problem can be formulated in the 
samplet domain, within a multi-kernel setting, or in a combined framework 
involving both, with the corresponding compressed operators.
This continuation progressively relaxes the sparsity penalty, stabilizes the
early iterations, and improves both robustness and convergence speed of the
TR-SSN method. This procedure has been investigated earlier in 
\cite{BNS97}. We remark that this approach of successively increasing the
number of active coefficients is necessary to make the use of the
block online SVD from Algorithm~\ref{alg:onlineSVD_matrix} computationally
feasible. Otherwise, one would have to resort to an iterative solver for
the computation of the Newton update, whose convergence might be critical
due to the ill-conditioned kernel matrices involved.

\subsection{Validation of tree-adaptive subsampling}
\label{sec:adaptive_sub_validation}

We validate the tree-adaptive subsampling strategy of
Section~\ref{sec:adaptive_ls} in isolation. 
The goal is to assess its ability to select a reduced set of representative
centers that capture the essential structure of the data, independently of 
any sparsity-promoting regularization. 
We compare the two energy measures based on the $\Xcal'$-norm and the 
$\Hcal'$-norm, which guide the adaptive selection and lead to different 
subsampling patterns. 

\begin{figure}[htb]
  \centering
  \begin{tikzpicture}
    \node[anchor=south west, inner sep=0] (img) at (0,0)
      {\includegraphics[width=0.42\textwidth]{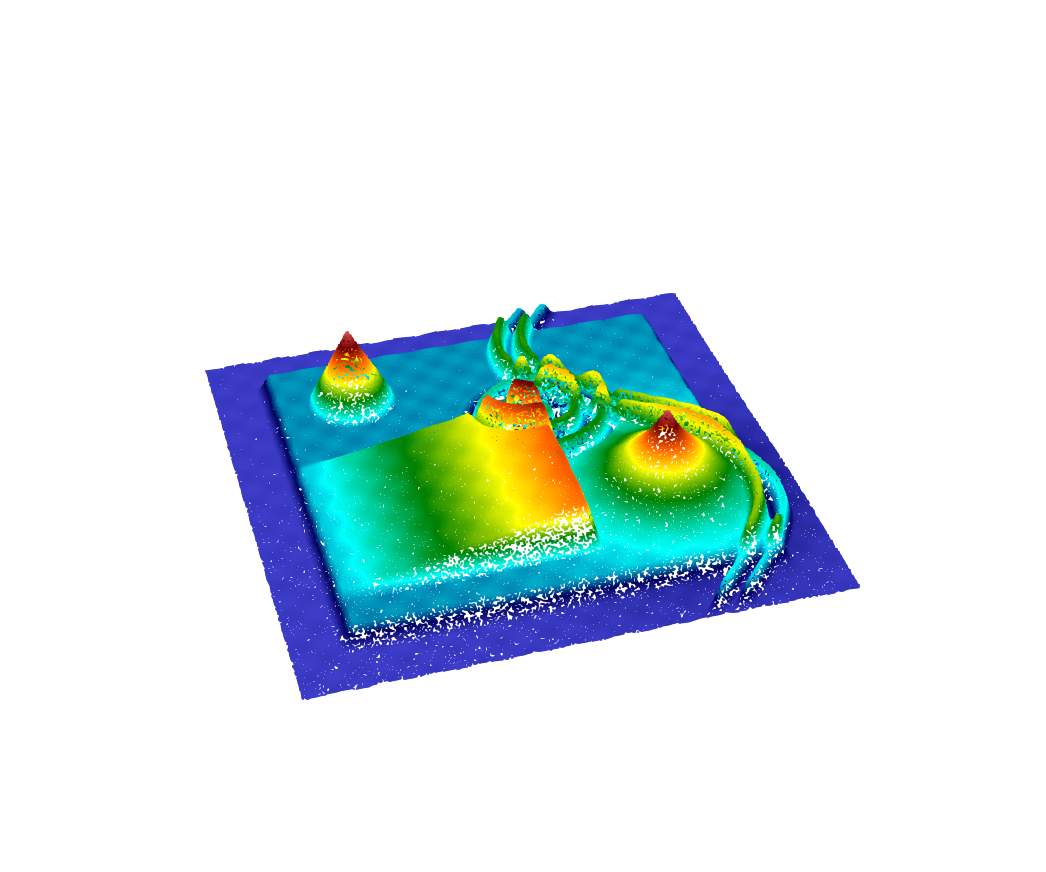}};
    \begin{scope}[
        x={(img.south east)},
        y={(img.north west)},
        font=\small,
        every path/.style={-{Stealth[length=5pt]}, thick},
        every node/.style={fill=white, fill opacity=0.75,
                           text opacity=1, inner sep=2pt,
                           rounded corners=2pt}
    ]
      \node (cone) at (0.05, 0.70) {cone peak};
      \draw (cone) -- (0.30, 0.62);
      \node (ridge) at (0.72, 0.85) {oscillatory wedge};
      \draw (ridge) -- (0.52, 0.57);
      \node (cusp) at (0.75, 0.20) {exponential cusp};
      \draw (cusp) -- (0.64, 0.49);
      \node (wedge) at (0.95, 0.65) {curved ridge};
      \draw (wedge) -- (0.65, 0.55);
      \node (grad) at (0.00, 0.35) {linear gradient};
      \draw (grad) -- (0.40, 0.45);
      \node (block) at (0.28, 0.85) {constant block};
      \draw (block) -- (0.43, 0.57);
    \end{scope}
  \end{tikzpicture}
  \caption{\textbf{Test 1.} Composite heterogeneous test function
$h$ defined on $[0,6]^2$.}\label{fig:ling_function}
\end{figure}
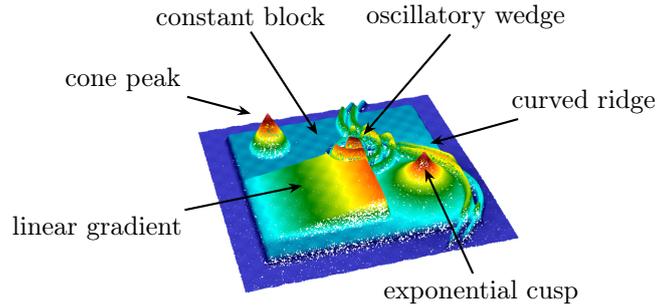

\paragraph{Test 1: Heterogeneous multiscale function}
We consider a very heterogeneous composite test function
$h\colon [0,6]^2 \to \mathbb{R}$, characterized by multiscale features
and localized irregularities, shown in Figure~\ref{fig:ling_function}.
The function is sampled on a large set of $N = 3\cdot10^5$ {uniformly distributed random}
data sites $X\subset [0,6]^2$, yielding the data vector
$\bm h = [h(\bm x_i)]_{i=1}^N$.
To exploit the multiscale structure of the data, we construct the
samplet tree associated with $X$ and compute the samplet coefficients
$\bm h^\Sigma = \bm T \bm h$.
An adaptive tree search is then performed with threshold {$t=\varepsilon^2\|{\bs h}^\Sigma\|_2^2$ in the $\Xcal'$-norm or
$t=\varepsilon^2({\bs h}^\Sigma)^\intercal{\bs D}_w{\bs h}^\Sigma$ in the
$\Hcal'$-norm,} see Section~\ref{sec:adaptive_ls}.
This procedure retains the most
significant coefficients of $\bm h^\Sigma$ and yields a reduced set of
representative centers $X_t\subset X$.
Given ${X_t}$, we construct a kernel approximation using the
exponential kernel
\begin{equation}\label{eq:kernel_exp}
    \kernel(\bm x, \bm y) = 
    \exp\bigg(-\frac{\|\bm x - \bm y\|_2}{\ell}\bigg),
\end{equation}
with lengthscale $\ell$
{chosen proportional to the fill distance $h_{X}$ of $X$, i.e., {$\ell = 10h_{X}$}. The same kernel defines the
$\Hcal'$-norm used for the selection.}
This choice balances accuracy and conditioning. {Smaller values of $\ell$ give a
more localized kernel and better conditioned matrices, while larger values
approximate smooth features better but worsen the conditioning.}
This leads to the kernel matrix
{${\bs K}_{X,X_t}\in\mathbb{R}^{N\times M}$} with $M = |X_t|$ and entries
$\kernel(\bm x_i, \widetilde{\bm x}_j)$,
$\bm x_i \in X$, $\widetilde{\bm x}_j \in X_t$.
The resulting least-squares problem is formulated in samplet coordinates,
and is solved by a QR decomposition.
This setting allows us to isolate
the effect of the adaptive subsampling strategy, without introducing
any sparsity-promoting regularization.
The experiment is repeated for several values of the threshold
$\varepsilon^2$ and for both subsampling criteria, enabling a direct
comparison between the $\Xcal'$-norm and $\Hcal'$-norm selections.

\begin{figure}[htb]
  \centering
  \includegraphics[width=0.24\linewidth]{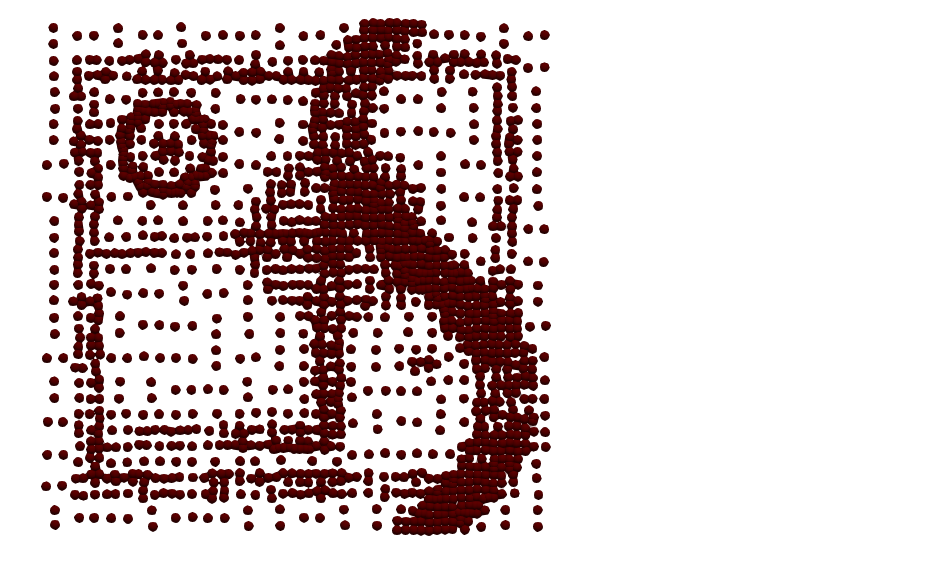}
  \includegraphics[width=0.24\linewidth]{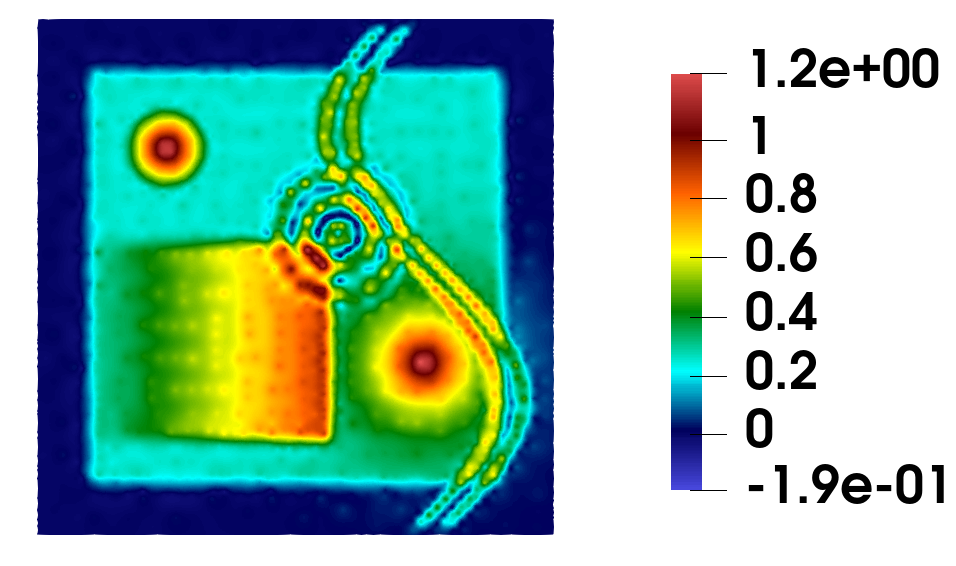}
  \qquad 
  \includegraphics[width=0.24\linewidth]{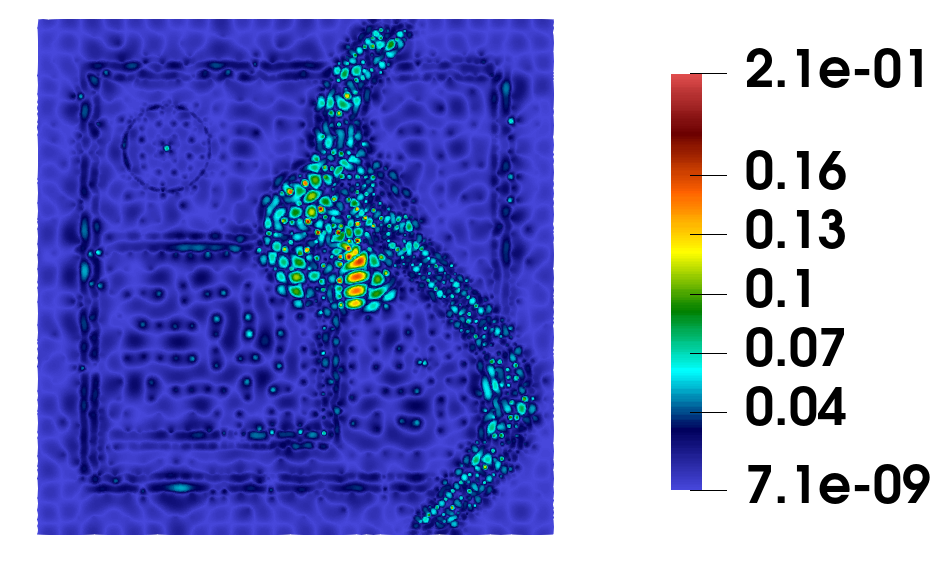}
  \\[1.0em]
  \includegraphics[width=0.24\linewidth]{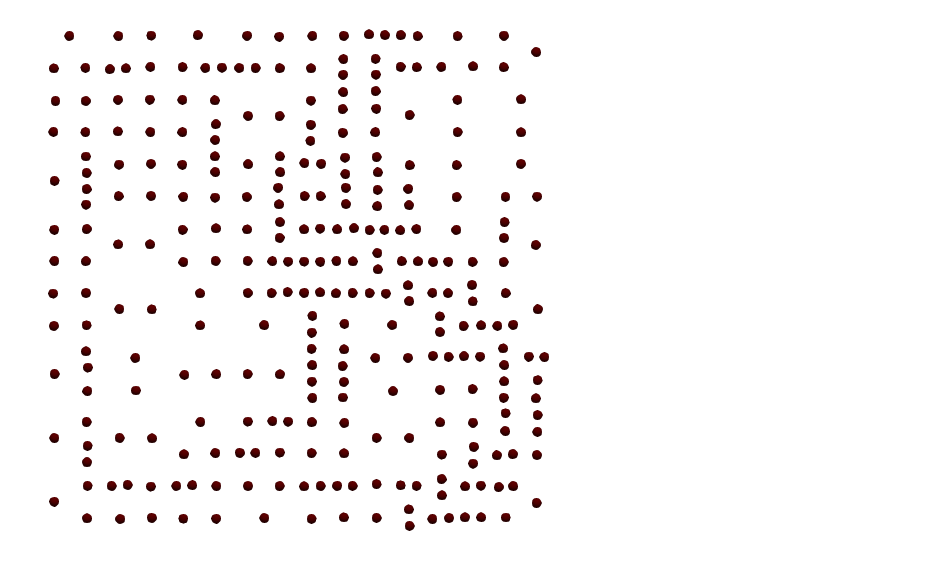}
  \includegraphics[width=0.24\linewidth]{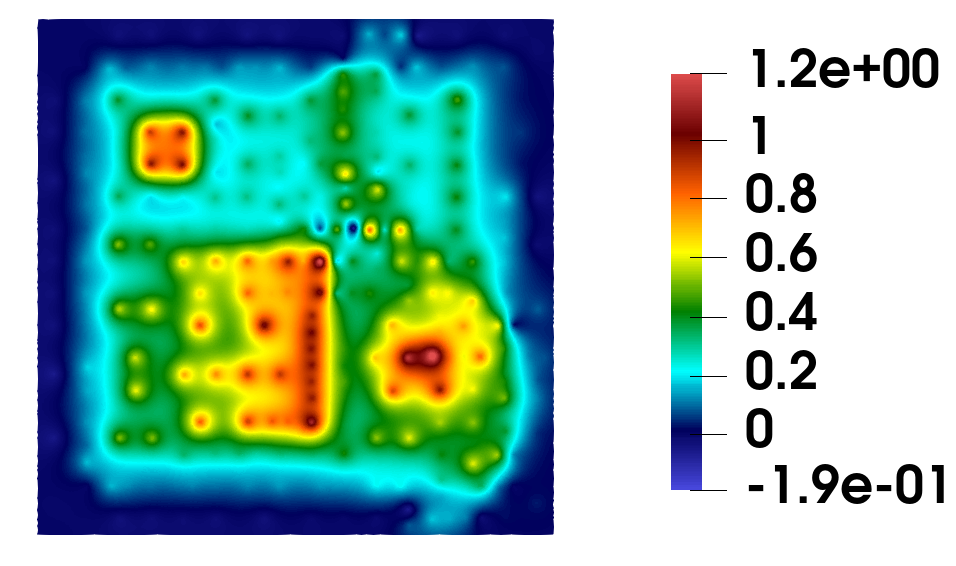}
  \qquad 
  \includegraphics[width=0.24\linewidth]{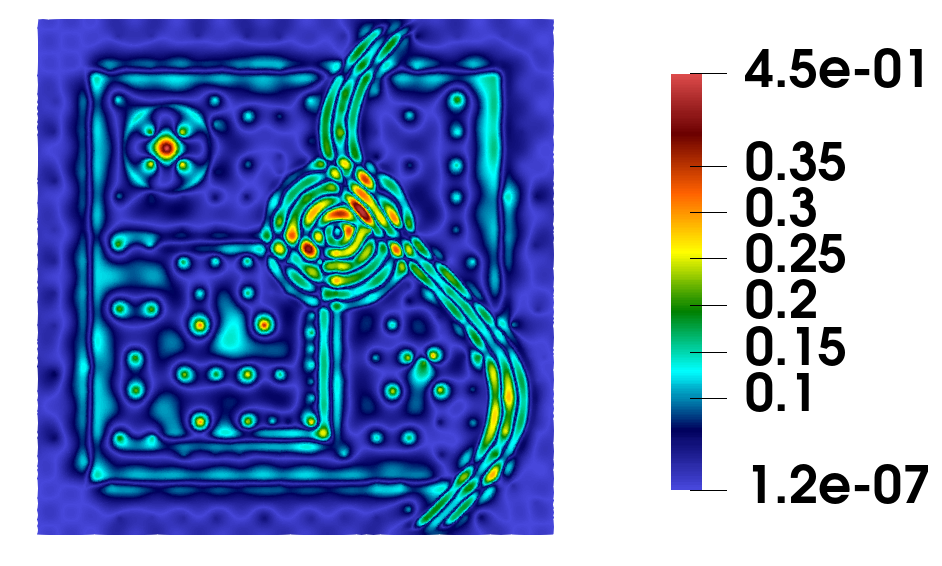}
\caption{\textbf{Test 1.} Reconstruction results for $\varepsilon^2 = 10^{-7}$.
  The two rows correspond to the $\Xcal'$-norm (top) and the $\Hcal'$-norm (bottom).
  From left to right, the columns display the adaptively selected points $X_t$,
  the reconstructed function and the pointwise absolute error.}
  \label{fig:adaptive_sub_reconstruction}
\end{figure}

A summary of the selected subsets and
the resulting reconstruction errors for different threshold values 
$\varepsilon^2$ is given in Table~\ref{tab:adaptive_sub}.
The reconstructed function is evaluated on a set of 
$N_{\mathrm{eval}} = 10^6$ uniformly distributed test points, and the
relative $\ell^2$-error is used to quantify the approximation quality.
Consistent with the discussion in Section~\ref{sec:adaptive_ls}, the
$\Hcal'$-norm criterion produces significantly sparser subsets, 
as it concentrates points in regions of high kernel
energy, whereas the $\Xcal'$-norm yields a more conservative selection. 
{For the chosen kernel, the ratio in \eqref{eq:relerrH} equals $5.7$, i.e., the relative
error in the $\Hcal'$-norm is bounded by $\sqrt{5.7}\,\varepsilon\approx2.4\,\varepsilon$.}
Figure~\ref{fig:adaptive_sub_reconstruction} illustrates the selected sets
$X_t$, the corresponding reconstructions, and the pointwise
absolute errors for $\varepsilon^2 = 10^{-7}$, comparing the
two subsampling criteria.

\begin{table}[htb]
  \centering
  \renewcommand{\arraystretch}{1.2}
  \setlength{\tabcolsep}{8pt}
  {\begin{tabular}{c c cc cc}
    \toprule
    & & \multicolumn{2}{c}{$\Xcal'$-norm}
      & \multicolumn{2}{c}{$\Hcal'$-norm} \\
    \cmidrule(lr){3-4} \cmidrule(lr){5-6}
    & $\varepsilon^2$ & $|X_t|$ & $e_2$ & $|X_t|$ & $e_2$ \\
    \midrule
    Test 1-A & $10^{-4}$  & 265  & $2.0\cdot10^{-1}$ & 53   & $5.2\cdot10^{-1}$ \\
    Test 1-B & $10^{-7}$  & 1696 & $4.7\cdot10^{-2}$ & 289  & $1.9\cdot10^{-1}$ \\
    Test 1-C & $10^{-10}$ & 4603 & $1.2\cdot10^{-2}$ & 1028 & $9.7\cdot10^{-2}$ \\
    \bottomrule
  \end{tabular}}
  \caption{\textbf{Test 1.} Tree-adaptive subsampling validation comparing
    the $\Xcal'$-norm and $\Hcal'$-norm criteria in terms of subset size
    $|X_t|$ and relative $\ell^2$ reconstruction error $e_2$, for varying
    threshold $\varepsilon^2$.}
  \label{tab:adaptive_sub}
\end{table}

\subsection{Validation of the full TR-SSN framework}\label{sec:TR_validation}

We now evaluate the full framework, combining tree-adaptive subsampling, 
samplet-based kernel compression,  and $\ell^1$-regularized regression solved 
via the TR-SSN method, see Section~\ref{sec:SSN}. 
The approach is tested on problems of increasing complexity:
a sparse sum of Gaussian functions in two dimensions (Test 2), 
the heterogeneous multiscale function from Figure \ref{fig:ling_function} (Test 3), 
and a large-scale 3D reflectance problem on the Stanford Bunny model \footnote{%
\url{https://graphics.stanford.edu/data/3Dscanrep/}} (Test 4). 
Varying parameters across tests are summarized in 
Table~\ref{tab:all_tests_summary}.

\begin{table}[htb]
\centering
\footnotesize
\setlength{\tabcolsep}{8pt}
\renewcommand{\arraystretch}{1.2}
\begin{tabular}{@{}l|ccccc@{}}
\toprule
  & T2-A & T2-B & T2-C & T3 & T4 \\
\midrule
\multicolumn{6}{c}{\textsc{General settings}} \\
\midrule
$d$
  & 2 & 2 & 2 & 2 & 3 \\
$N$
  & $10^6$ & $10^6$ & $10^6$
  & $10^6$ & {$1.2\cdot10^6$}  \\
\addlinespace[4pt]
\multicolumn{6}{c}{\textsc{Adaptive subsampling}} \\
\midrule
{Norm}
  & {$\Hcal'$, $\ell=2\rho_X$} & {$\Hcal'$, $\ell=2\rho_X$} & {$\Hcal'$, $\ell=2\rho_X$}
  & {$\Hcal'$, $\ell=4\rho_X$} & {$\Hcal'$, $\ell=2\rho_X$} \\
$\varepsilon^2$
  & $10^{-12}$ & $10^{-12}$ & $10^{-12}$
  & $10^{-8}$ & $10^{-5}$ \\
$|X_t|$
  & {856} & {856} & {856}
  & {2541} & {4857} \\
\addlinespace[4pt]
\multicolumn{6}{c}{\textsc{Kernel model}} \\
\midrule
Kernel type
  & Gauss. & Gauss. & Gauss.
  & exp. & exp. \\
$L$
  & 1 & 1 & 4
  & 5 & 3 \\
$\ell_j \ (j = 1,\dots,L)$
  & $h_{X_t}$ & $h_{X_t}$ & $\rho_{X_t}
\left(\frac{2 h_{X_t}}{\rho_{X_t}}\right)^{\frac{j-1}{L-1}}$
  & $2\rho_{X_t}
\left(\frac{2 h_{X_t}}{2\rho_{X_t}}\right)^{\frac{j-1}{L-1}}$ & $5\rho_{X_t}
\left(\frac{2 h_{X_t}}{5\rho_{X_t}}\right)^{\frac{j-1}{L-1}}$ \\
\addlinespace[4pt]
\multicolumn{6}{c}{\textsc{TR-SSN}} \\
\midrule
$\mathtt{tol}$
  & $10^{-5}$ & $10^{-5}$ & $10^{-5}$
  & $10^{-6}$ & $10^{-5}$ \\
$\mathtt{maxit}$
  & 100 & 100 & 100
  & 50 & 50 \\
\addlinespace[4pt]
\multicolumn{6}{c}{\textsc{Continuation strategy}} \\
\midrule
$r_0$
  & 10 & 10 & 10
  & 10 & 10 \\
$r_{\min}$
  & $10^{-5}$ & $10^{-5}$ & $10^{-5}$
  & {$10^{-4}$} & $10^{-4}$ \\
$\gamma$
  & 0.7 & 0.7 & 0.7
  & 0.7 & 0.75 \\
\bottomrule
\end{tabular}
\caption{
Summary of all parameters used across Tests~2--4.
Test~2 compares three formulations: dense single-kernel (A), samplet-compressed
single-kernel (B), and samplet-compressed multi-kernel (C).
Here, 
{$\rho_X$ and $\rho_{X_t}$ denote the separation radii of $X$ and
$X_t$}, i.e., the minimum distance between distinct {points}, while
$h_{X_t}$ is the fill distance, measuring how well $X_t$ covers the domain.}
\label{tab:all_tests_summary}
\end{table}
\paragraph{Test 2: Multiscale kernel regression benchmark}
In this experiment, we assess the performance of the full TR-SSN framework
on a multiscale regression problem in two dimensions.
We consider a set $X$ of $N = 10^6$ uniformly distributed random points,
where we sample the function $h \colon [0,1]^2 \to \mathbb{R}$ defined as a
superposition of four Gaussian components
\[
h(\bm x)
= \sum_{k=1}^{4} a_k \exp\!\left(
-\frac{\|\bm x - \bm c_k\|_2^2}{2\sigma_k^2}
\right), \qquad \bm x \in [0,1]^2,
\]
with centers $\bm c_k$, amplitudes $a_k$, and widths $\sigma_k$ given by
\[
\begin{array}{lll}
\bm c_1 = (0.20,\, 0.80)^\top, & a_1 = 0.50, & \sigma_1 = 0.03, \\
\bm c_2 = (0.75,\, 0.25)^\top, & a_2 = 0.60, & \sigma_2 = 0.06, \\
\bm c_3 = (0.65,\, 0.65)^\top, & a_3 = 0.40, & \sigma_3 = 0.13, \\
\bm c_4 = (0.30,\, 0.30)^\top, & a_4 = 0.55, & \sigma_4 = 0.25.
\end{array}
\]
Starting from the data vector $\bm h = [h(\bm x_i)]_{i=1}^N$, we first apply
tree-adaptive subsampling based on the $\Hcal'$-norm with an exponential kernel of lengthscale twice
the separation radius of $X$. {The indicator in
\eqref{eq:relerrH} amounts to $19.4$.
This yields} a reduced set of
representative centers $X_t \subset X$ of cardinality {$|X_t| = 856$}.
We then construct kernel-based approximations using Gaussian kernels
$\kernel(\bm x, \bm y)
= \exp\left(-{\|\bm x - \bm y\|_2^2}/({2\, \ell^2})\right)$, combined,
{in Tests~2-B and~2-C,} with samplet compression as described in
Section~\ref{sec:samplets}, and solve the resulting $\ell^1$-regularized
problems using the TR-SSN method in Algorithm~\ref{alg:trustSSN}.
We compare three different formulations.
In the first case (Test 2-A), we consider a dense single-kernel model and solve
the standard $\ell^1$-regularized least-squares problem~\eqref{eq:SSLASSO}.
In the second case (Test 2-B), we employ a samplet-compressed single-kernel
formulation and solve the problem in samplet coordinates~\eqref{eq:MSLASSO}.
In the third case (Test 2-C), we extend the compressed formulation to a
multi-kernel setting, using a dictionary of 4 kernels with lengthscales from the separation radius to twice
the fill distance of $X_t$.

\begin{figure}[htb]
    \centering
        \begin{subfigure}[b]{0.33\textwidth}
        \centering
        \includegraphics[width=0.75\textwidth]{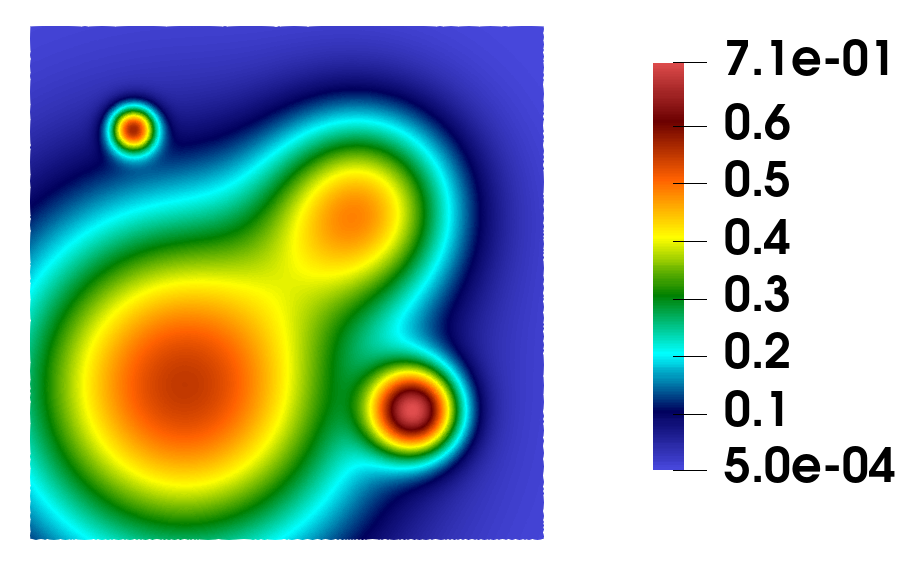}
    \end{subfigure}
    \begin{subfigure}[b]{0.33\textwidth}
        \centering
        \includegraphics[width=0.75\textwidth]{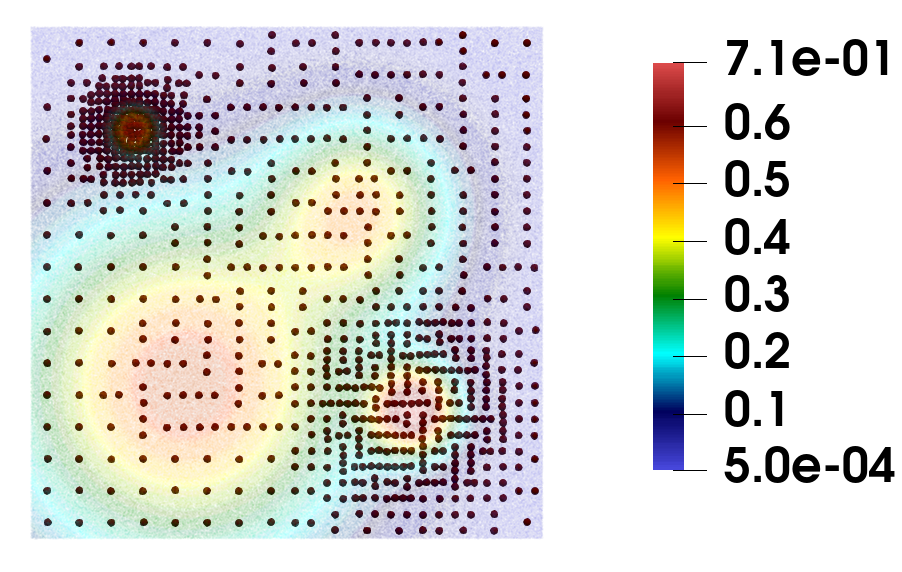}
    \end{subfigure} \\
    \vspace{0.4 cm}
    
    \begin{subfigure}[b]{0.33\textwidth}
        \centering
        \includegraphics[width=0.75\textwidth]{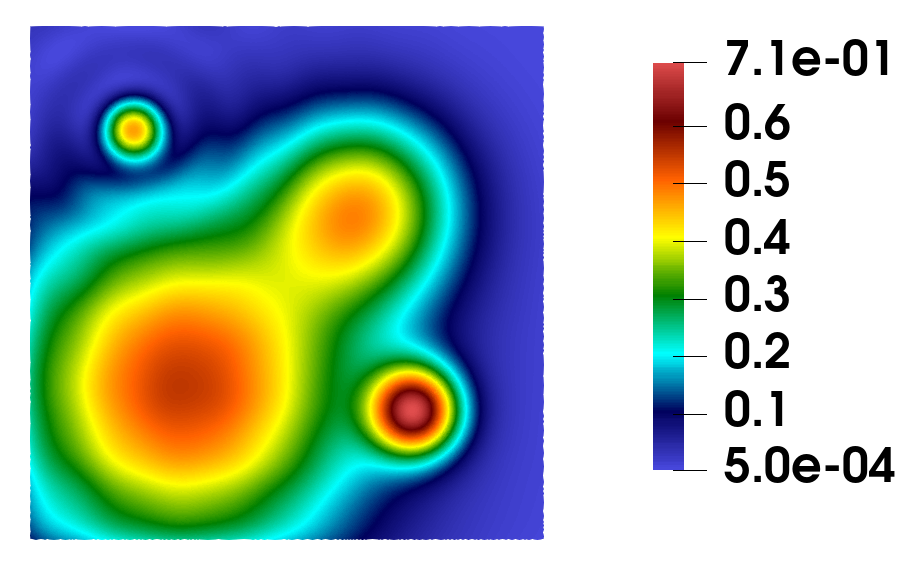}
    \end{subfigure}
    \hspace{-0.5 cm}
    \begin{subfigure}[b]{0.33\textwidth}
        \centering
        \includegraphics[width=0.75\textwidth]{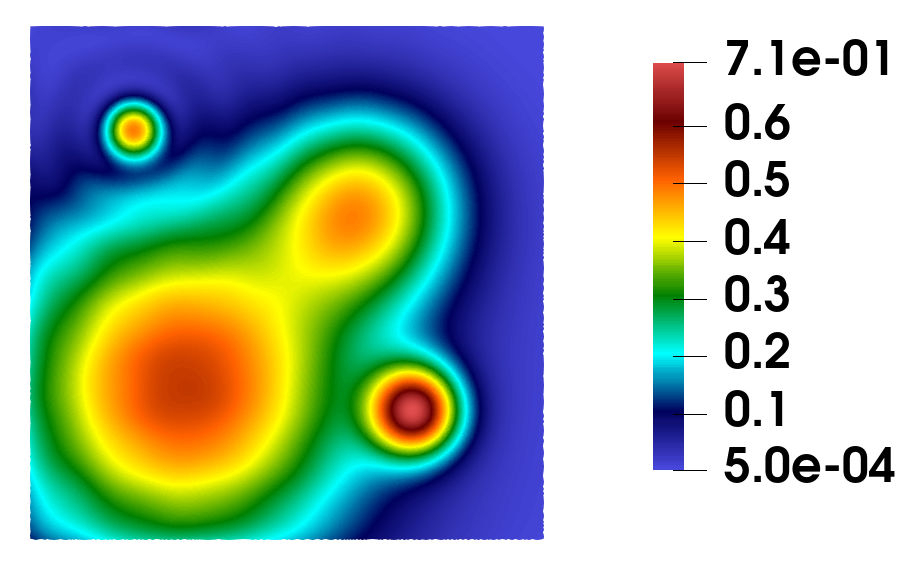}
    \end{subfigure}
        \hspace{-0.5 cm}
  \begin{subfigure}[b]{0.33\textwidth}
        \centering
        \includegraphics[width=0.75\textwidth]{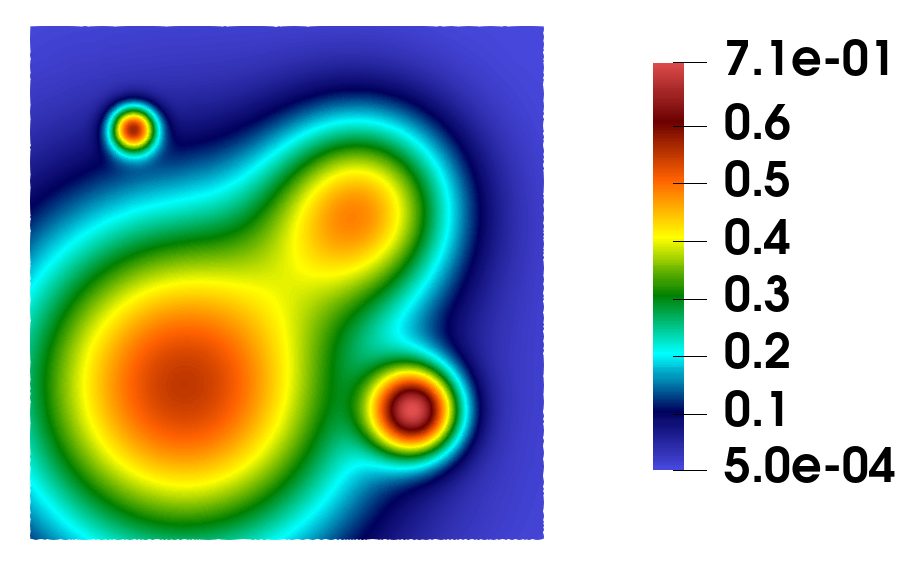}
    \end{subfigure}
    \vskip\baselineskip  
    \begin{subfigure}[b]{0.33\textwidth}
        \centering
        \includegraphics[width=0.75\textwidth]{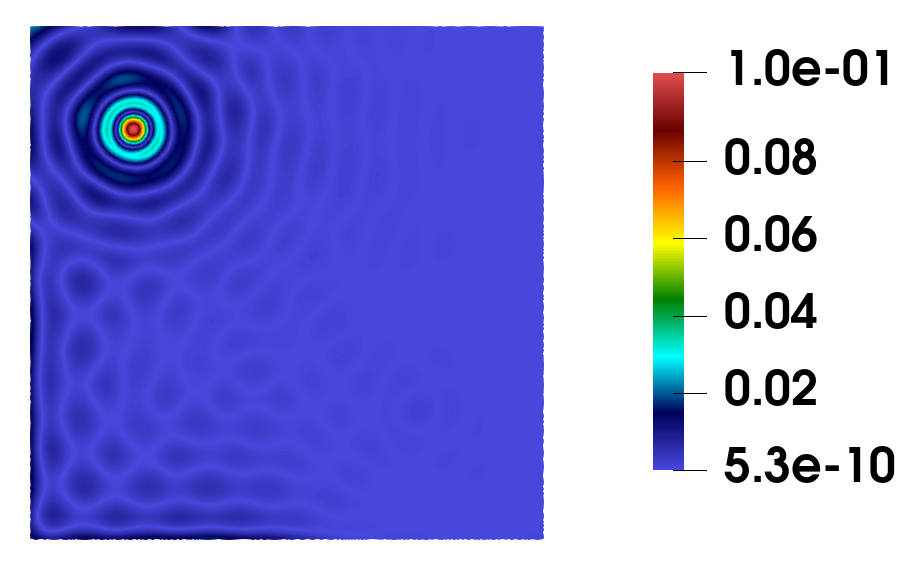}
    \end{subfigure}
    \hspace{-0.5 cm}
    \begin{subfigure}[b]{0.33\textwidth}
        \centering
        \includegraphics[width=0.75\textwidth]{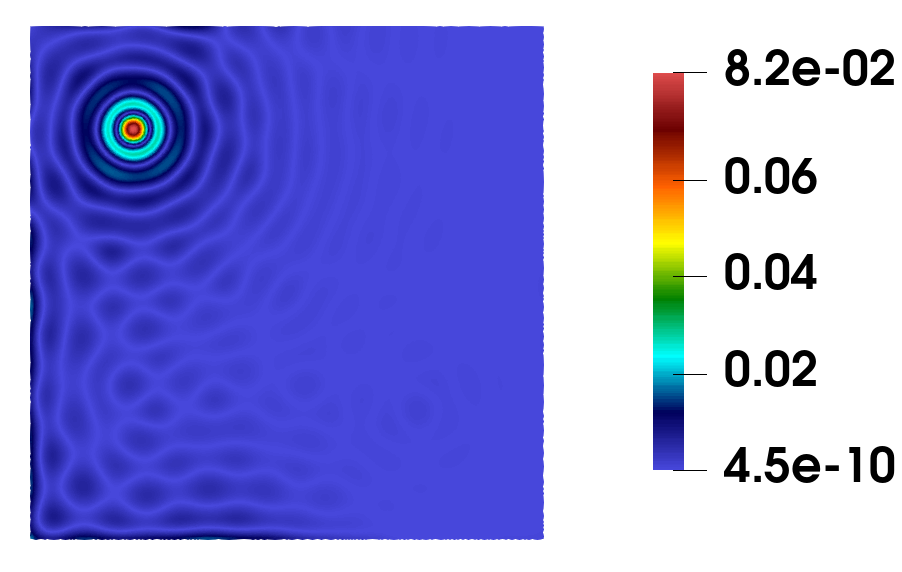}
    \end{subfigure}
        \hspace{-0.5 cm}
    \begin{subfigure}[b]{0.33\textwidth}
        \centering
        \includegraphics[width=0.75\textwidth]{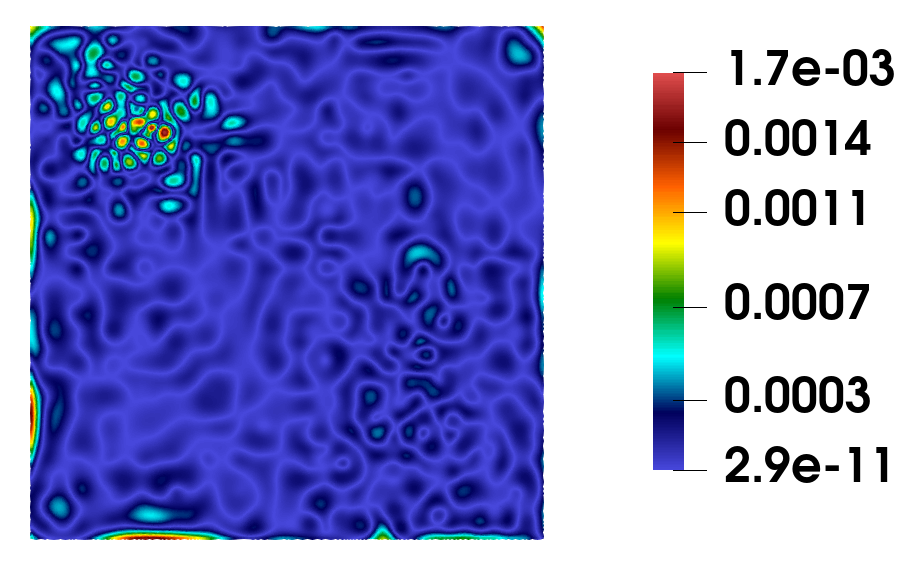}
    \end{subfigure}
\caption{\textbf{Test 2.} Visualization of the multiscale kernel regression
  benchmark. The first row shows the reference function and the adaptive 
  samplet-selected points. 
The second row displays the reconstructed functions for the dense single-kernel 
(Test 2-A), samplet-compressed single-kernel (Test 2-B), and samplet-compressed 
multi-kernel models (Test 2-C). 
The third row reports the corresponding pointwise absolute errors for each 
reconstruction. }
    \label{fig:res_test2}
\end{figure}

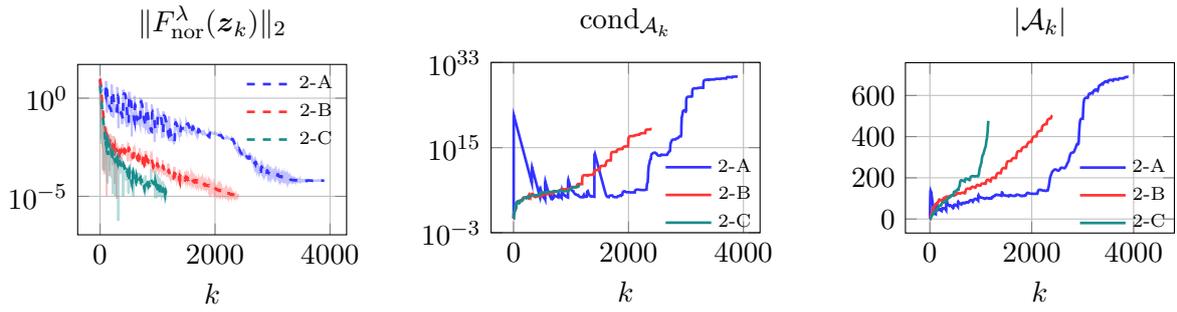
\begin{figure}[htb]
\centering
\begin{minipage}{0.32\textwidth}
\centering
\begin{tikzpicture}
\begin{semilogyaxis}[
    width=\textwidth,
    height=0.8\textwidth,
    xlabel={$k$},
    title={$\|F_{\text{nor}}^\lambda(\bs z_k)\|_2$},
    grid=both,
    grid style={line width=0.2pt, draw=gray!30},
    major grid style={line width=0.4pt, draw=gray!50},
    scaled ticks=false,
    tick label style={
        /pgf/number format/1000 sep={}
    },
    legend style={
    draw=none,
    fill=none,
    font=\tiny,
    at={(1.0,1.0)},anchor=north east
    }
]

\addplot[blue!60, line width=0.8pt, opacity=0.4, forget plot]
table[x=iter,y=residual]{./figures/rev_benchmark1_pgf.txt};
\addplot[blue!80, dashed, line width=1pt]
table[x=iter,y=residual_smooth]{./figures/rev_benchmark1_pgf.txt};
\addlegendentry{2-A}

\addplot[red!60, line width=0.8pt, opacity=0.4, forget plot]
table[x=iter,y=residual]{./figures/rev_benchmark2_pgf.txt};
\addplot[red!80, dashed, line width=1pt]
table[x=iter,y=residual_smooth]{./figures/rev_benchmark2_pgf.txt};
\addlegendentry{2-B}

\addplot[teal!70, line width=0.8pt, opacity=0.4, forget plot]
table[x=iter,y=residual]{./figures/rev_benchmark3_pgf.txt};
\addplot[teal!90, dashed, line width=1pt]
table[x=iter,y=residual_smooth]{./figures/rev_benchmark3_pgf.txt};
\addlegendentry{2-C}

\end{semilogyaxis}
\end{tikzpicture}
\end{minipage}
\hfill
\begin{minipage}{0.32\textwidth}
\centering
\begin{tikzpicture}
\begin{semilogyaxis}[
    width=\textwidth,
    height=0.8\textwidth,
    xlabel={$k$},
    title={$\operatorname{cond}_{\Act_k}$},
    grid=both,
    grid style={line width=0.2pt, draw=gray!30},
    major grid style={line width=0.4pt, draw=gray!50},
    scaled ticks=false,
    tick label style={
        /pgf/number format/1000 sep={}
    },
    legend style={
    draw=none,
    fill=none,
    font=\tiny,
    at={(0.5,1.0)},anchor=north east
    }
]

\addplot[blue!80, line width=1pt]
table[x=iter,y=cond_smooth]{./figures/rev_benchmark1_pgf.txt};
\addlegendentry{2-A}

\addplot[red!80, line width=1pt]
table[x=iter,y=cond_smooth]{./figures/rev_benchmark2_pgf.txt};
\addlegendentry{2-B}

\addplot[teal!90, line width=1pt]
table[x=iter,y=cond_smooth]{./figures/rev_benchmark3_pgf.txt};
\addlegendentry{2-C}

\end{semilogyaxis}
\end{tikzpicture}
\end{minipage}
\hfill
\begin{minipage}{0.32\textwidth}
\centering
\begin{tikzpicture}
\begin{axis}[
    width=\textwidth,
    height=0.8\textwidth,
    xlabel={$k$},
    title={$|\Act_k|$},
    grid=both,
    grid style={line width=0.2pt, draw=gray!30},
    major grid style={line width=0.4pt, draw=gray!50},
    scaled ticks=false,
    tick label style={
        /pgf/number format/1000 sep={}
    },
    xmax = 4300,
    legend style={
    draw=none,
    fill=none,
    font=\tiny,
    at={(0.45,1.0)},anchor=north east
    }
]

\addplot[blue!80, line width=1pt]
table[x=iter,y=nactive_smooth]{./figures/rev_benchmark1_pgf.txt};
\addlegendentry{2-A}

\addplot[red!80, line width=1pt]
table[x=iter,y=nactive_smooth]{./figures/rev_benchmark2_pgf.txt};
\addlegendentry{2-B}

\addplot[teal!90, line width=1pt]
table[x=iter,y=nactive_smooth]{./figures/rev_benchmark3_pgf.txt};
\addlegendentry{2-C}

\end{axis}
\end{tikzpicture}
\end{minipage}

\caption{\textbf{Test 2.} TR-SSN iteration metrics. Faint lines show the 
raw data, while dashed lines represent a smoothed trend obtained via an 
exponential moving average, which reduces short-term fluctuations.}
\label{fig:test2metrics}
\end{figure}

The resulting reconstructions, selected centers, and pointwise absolute errors
are shown in Figure~\ref{fig:res_test2}, while quantitative results in terms of
sparsity and relative $\ell^2$ error {measured on all data sites $X$}
are reported in Table~\ref{tab:test2_sparsity}. Figure~\ref{fig:test2metrics}
shows representative TR-SSN iteration metrics for the three configurations.
{The single-kernel models of Tests 2-A and 2-B provide} a reasonable
approximation but lack flexibility in capturing the different scales of the
target function. {Moreover, for the dense formulation, TR-SSN reaches the
maximum number of iterations on every continuation level. On the final levels,
its residual stagnates at about $10^{-3}$, due to the growth of the condition
number of the matrices on the active set beyond $10^{29}$.}
{Compared to the dense formulation,} the samplet-based formulations
improve both accuracy and sparsity, with the multi-kernel approach achieving
the smallest error while maintaining a compact representation of the solution.

\begin{table}[htb]
\centering
\footnotesize
\setlength{\tabcolsep}{12pt}
{\begin{tabular}{@{}lcc@{}}
\toprule
 & $\|\bm\alpha_j\|_0/|X_t|$ & $e_2$ \\
\midrule
\multicolumn{3}{c}{\textbf{Test 2-A}} \\[4pt]
$\bm K$ ($\ell = 6.53\cdot10^{-2}$) & $840/856$ & $1.87\cdot10^{-2}$ \\[8pt]
\multicolumn{3}{c}{\textbf{Test 2-B}} \\[4pt]
$\bm K$ ($\ell = 6.53\cdot10^{-2}$) & $612/856$ & $1.48\cdot10^{-2}$ \\[8pt]
\multicolumn{3}{c}{\textbf{Test 2-C}} \\[4pt]
$\bm K_1$ ($\ell_1 = 6.13\cdot10^{-3}$) & $5 / 856$   & \\[4pt]
$\bm K_2$ ($\ell_2 = 1.70\cdot10^{-2}$) & $227 / 856$ & \multirow{3}{*}{$5.80\cdot10^{-4}$}\\[4pt]
$\bm K_3$ ($\ell_3 = 4.71\cdot10^{-2}$) & $206 / 856$ & \\[4pt]
$\bm K_4$ ($\ell_4 = 1.31\cdot10^{-1}$) & $35 / 856$  & \\[4pt]
\bottomrule
\end{tabular}}
\caption{\textbf{Test 2.} Sparsity of the coefficient vectors and relative
$\ell^2$ reconstruction error for the three test cases.}
\label{tab:test2_sparsity}
\end{table}

\paragraph{Test 3: Sparse multi-kernel reconstruction of a heterogeneous 2D signal}
We consider again the heterogeneous multiscale function
$h\colon [0,6]^2 \to \mathbb{R}$ introduced in
Test 1, see Figure~\ref{fig:ling_function}, sampled on {$N=10^6$} {uniformly distributed random} data sites
$X$.
Starting from the data vector $\bm h = [h(\bm x_i)]_{i=1}^N$, we apply
tree-adaptive subsampling based on the $\Hcal'$-norm described in
Section~\ref{sec:adaptive_ls}, {with an exponential kernel of
lengthscale four times the separation radius of $X$. The indicator in
\eqref{eq:relerrH} amounts to $47.9$.}
It yields a reduced set of representative centers
$X_t \subset X$ {of cardinality {$|X_t| = 2541$}}.
On this reduced set, we construct a multi-kernel approximation using
a dictionary of 5 exponential kernels \eqref{eq:kernel_exp} {with lengthscales 
ranging from twice the separation radius to twice the fill distance of $X_t$,}
combined with the samplet-compressed representation introduced in
Section~\ref{sec:samplets}.
The resulting weighted $\ell^1$-regularized least-squares problem is
formulated in samplet coordinates and solved via the TR-SSN method, see
Algorithm~\ref{alg:trustSSN}.

\begin{figure}[htb]
    \centering
    \begin{subfigure}[b]{0.33\textwidth}
        \centering        \includegraphics[width=0.8\textwidth]{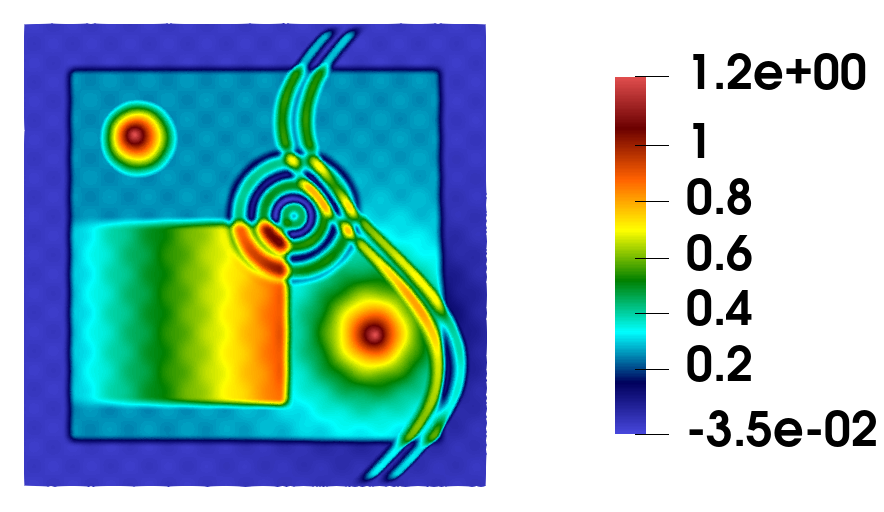}
    \end{subfigure}
    \hspace{0.5 cm}
  \begin{subfigure}[b]{0.33\textwidth}
        \centering
        \includegraphics[width=0.8\textwidth]{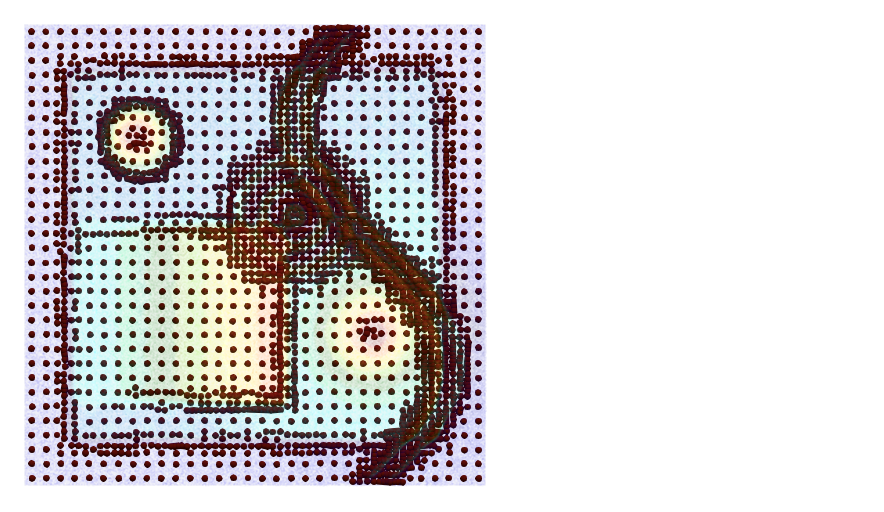}
    \end{subfigure}
    \vskip\baselineskip  
    \begin{subfigure}[b]{0.33\textwidth}
        \centering
        \includegraphics[width=0.8\textwidth]{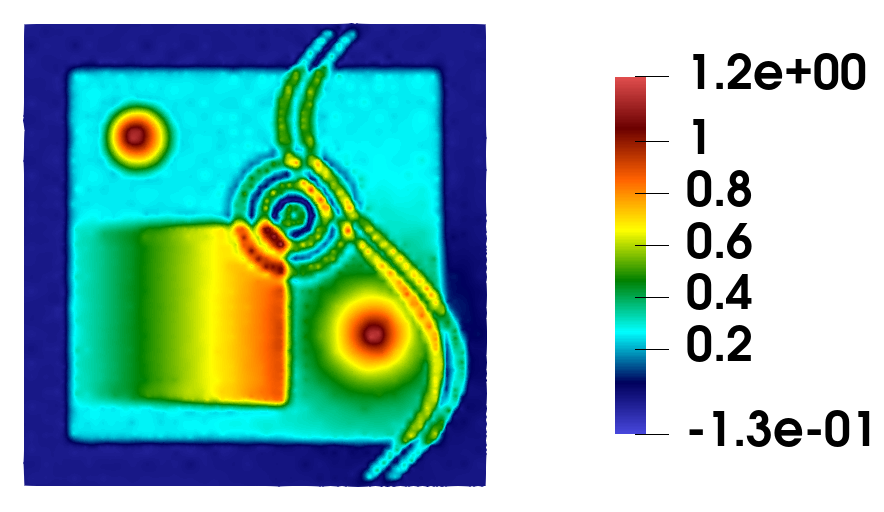}
    \end{subfigure}
    \hspace{-0.2 cm}
    \begin{subfigure}[b]{0.33\textwidth}
        \centering
        \includegraphics[width=0.8\textwidth]{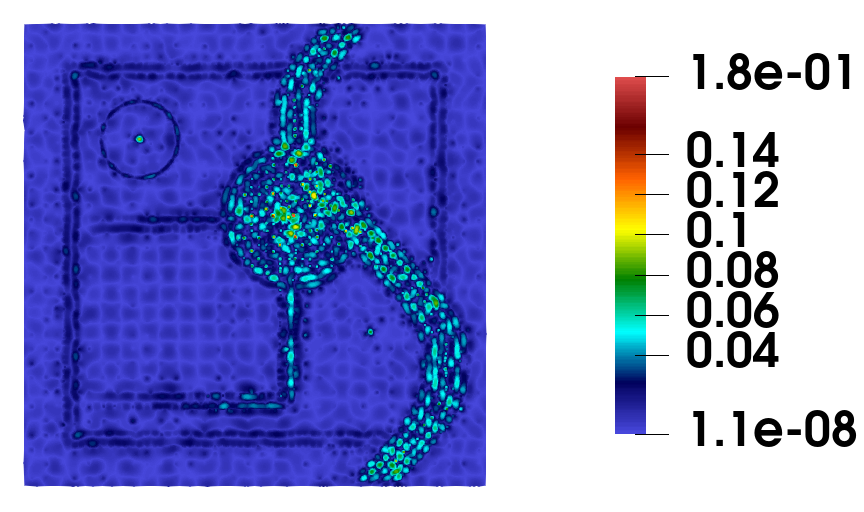}
    \end{subfigure}
    \hspace{-0.2 cm}
    \begin{subfigure}[b]{0.33\textwidth}
        \centering
        \includegraphics[width=0.8\textwidth]{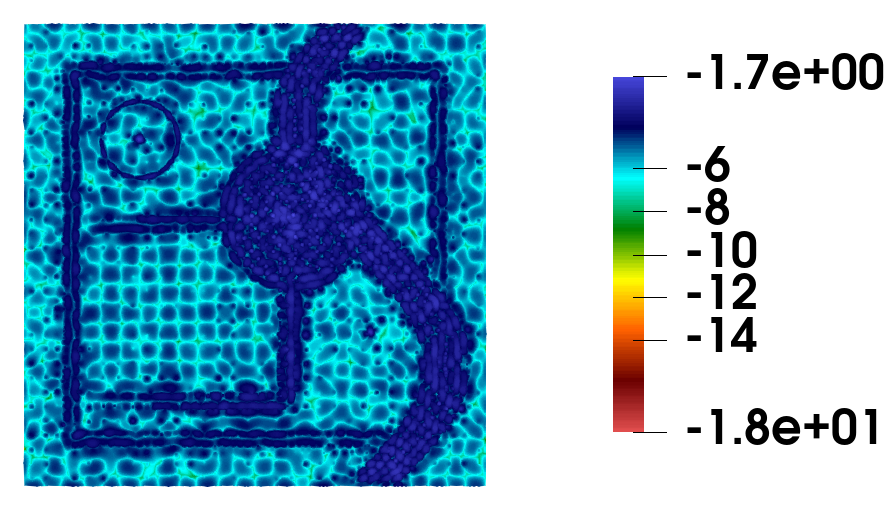}
    \end{subfigure}
    \begin{subfigure}[b]{0.33\textwidth}
        \centering
        \hspace{0.1 cm}\includegraphics[width=0.86\textwidth]{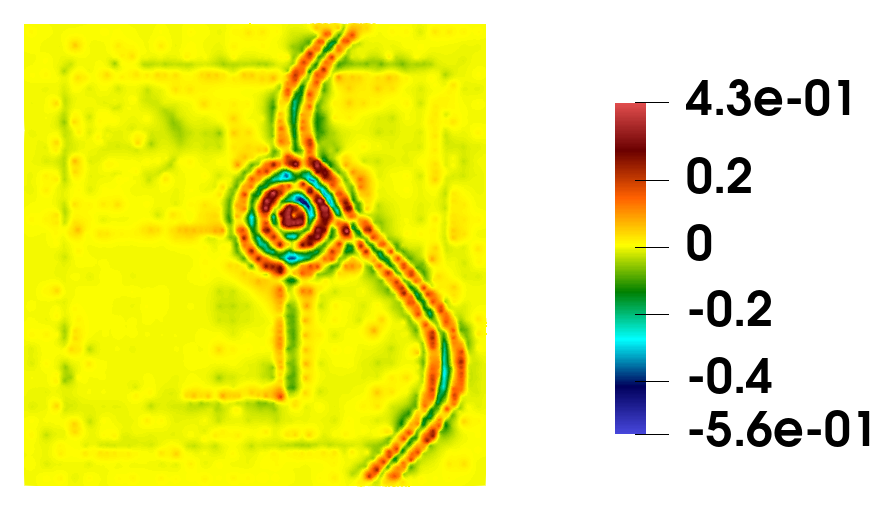}
    \end{subfigure}
    \hspace{-0.2 cm}
    \begin{subfigure}[b]{0.33\textwidth}
        \centering
        \includegraphics[width=0.8\textwidth]{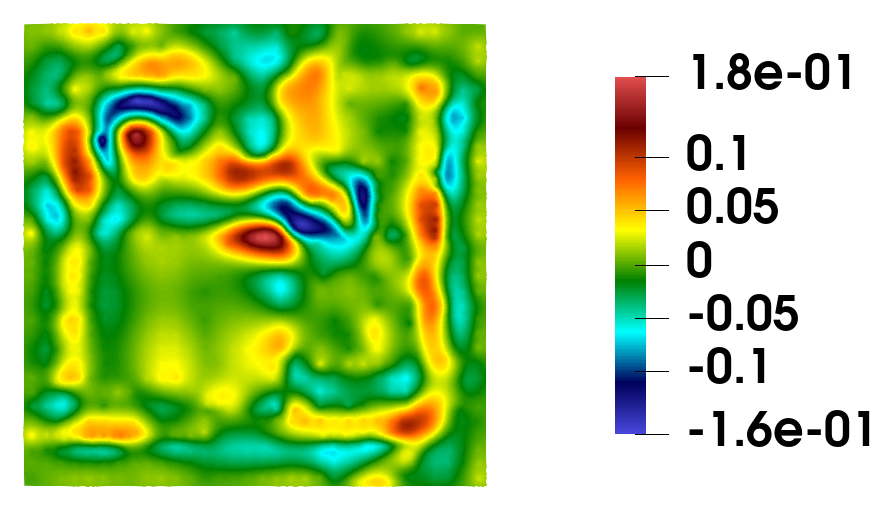}
    \end{subfigure}
    \hspace{-0.2 cm}
    \begin{subfigure}[b]{0.33\textwidth}
        \centering
        \includegraphics[width=0.8\textwidth]{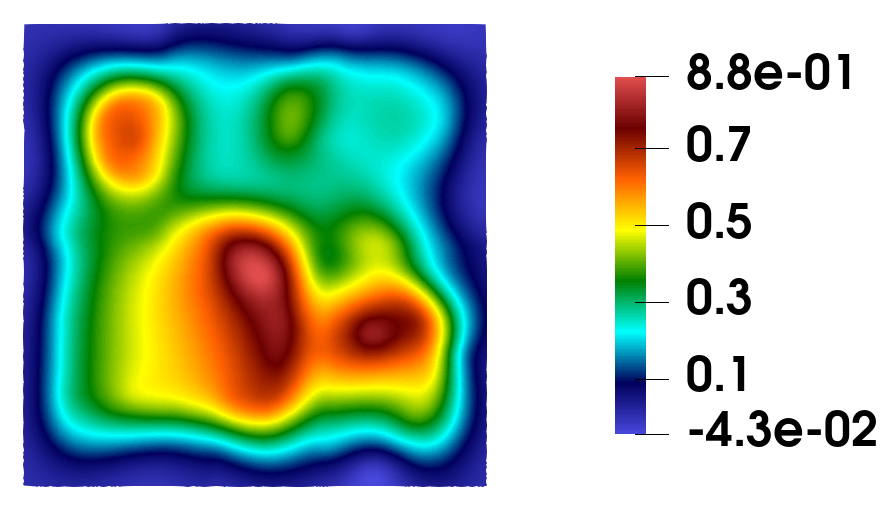}
    \end{subfigure}
\caption{\textbf{Test 3.} First row: reference solution and adaptive point set $X_t$. 
Second row: reconstructed solution, pointwise absolute error, and pointwise absolute error in logarithmic scale. 
Third row: contributions of the kernels $\bm K_1$, $\bm K_3$, and $\bm K_5$ in the multi-kernel expansion.}    \label{fig:res_test3}
\end{figure}

The selected centers, reconstructed solution, and corresponding errors
are shown in Figure~\ref{fig:res_test3}, while quantitative results in
terms of sparsity and reconstruction accuracy are reported in
Table~\ref{tab:test3_timing}.
The adaptive selection concentrates points in regions of high variability,
and the multi-kernel formulation effectively captures both fine and
coarse features of the signal, yielding an accurate reconstruction
with a sparse representation.

\begin{table}[htb]
\centering
\footnotesize
\setlength{\tabcolsep}{12pt}
{\begin{tabular}{@{}lcc@{}}
\toprule
 & $\|\bm\alpha_j\|_0 / |X_t|$ & $e_2$ \\
\midrule
$\bm K_1$ ($\ell_1 = 7.31\cdot10^{-2}$) & $1260 / 2541$ & \\[4pt]
$\bm K_2$ ($\ell_2 = 1.11\cdot10^{-1}$) & $574 / 2541$  & \multirow{3}{*}{$3.47\cdot10^{-2}$}\\[4pt]
$\bm K_3$ ($\ell_3 = 1.67\cdot10^{-1}$) & $259 / 2541$  & \\[4pt]
$\bm K_4$ ($\ell_4 = 2.53\cdot10^{-1}$) & $138 / 2541$  & \\[4pt]
$\bm K_5$ ($\ell_5 = 3.84\cdot10^{-1}$) & $84 / 2541$   & \\
\bottomrule
\end{tabular}}
\caption{\textbf{Test 3.} Sparsity of the coefficient vectors $\bm\alpha_1$--$\bm\alpha_5$
across the {five} kernel lengthscales and relative $\ell^2$ reconstruction error
for the composite function reconstruction ({$|X_t| = 2541$} centers).}
\label{tab:test3_timing}
\end{table}

\paragraph{Test 4: Large-scale 3D reflectance}
In this experiment, we consider a large-scale 3D surface problem
on the Stanford bunny model, represented by a point cloud with an associated scalar
reflectance field. The dataset $X$ contains the 3D coordinates of
{$N \approx 1.2 \cdot 10^6$} vertices, and the function
$h \colon X \to \mathbb{R}$ is generated by classical diffuse-specular reflectance,
see, e.g.,~\cite{phong1975illumination, foley1996computer}, i.e.,
{
\begin{equation}\label{eq:phong}
  h(\bm x)
  = 0.4 + \max\bigl(n(\bm x)\cdot\bm v_l,\,0\bigr)
  \bigg[0.5 + 0.1\exp\bigg(-\frac{\|\bm r(\bm x)-\bm v_o\|_2^2}{8}\bigg)\bigg],
\end{equation}
}
where $n(\bm x)$ denotes the unit outward normal at $\bm x$ and
$\bm r(\bm x) = 2\big(n(\bm x)\cdot\bm v_l\big)n(\bm x) - \bm v_l$
is the mirror reflection direction.
The light and view directions are fixed as
$\bm v_l = \tfrac{1}{\sqrt{2}}[1,\,0,\,1]^\top$ and
$\bm v_o = \tfrac{1}{\sqrt{6}}[1,\,1,\,2]^\top$.
Starting from the data vector $\bm h$, we perform tree-adaptive subsampling
{based on the $\Hcal'$-norm, see Section~\ref{sec:adaptive_ls}, with an
exponential kernel of lengthscale twice the separation radius of $X$. The
indicator in \eqref{eq:relerrH} amounts to $1.12$. This yields} a reduced set of representative
centers $X_t \subset X$ of cardinality {$|X_t| = 4857$}.
On this reduced set, we construct a multi-kernel approximation using
a small dictionary of 3 exponential kernels \eqref{eq:kernel_exp} with geometrically distributed
lengthscales from five times the separation radius to twice the fill distance of $X_t$, combined with samplet compression.
The resulting $\ell^1$-regularized problem~\eqref{eq:MSLASSO} is then solved via the TR-SSN method, see Algorithm~\ref{alg:trustSSN}.

\begin{figure}[htbp!]
    \centering
    \begin{subfigure}[b]{0.33\textwidth}
        \centering
        \includegraphics[width=0.85\textwidth]{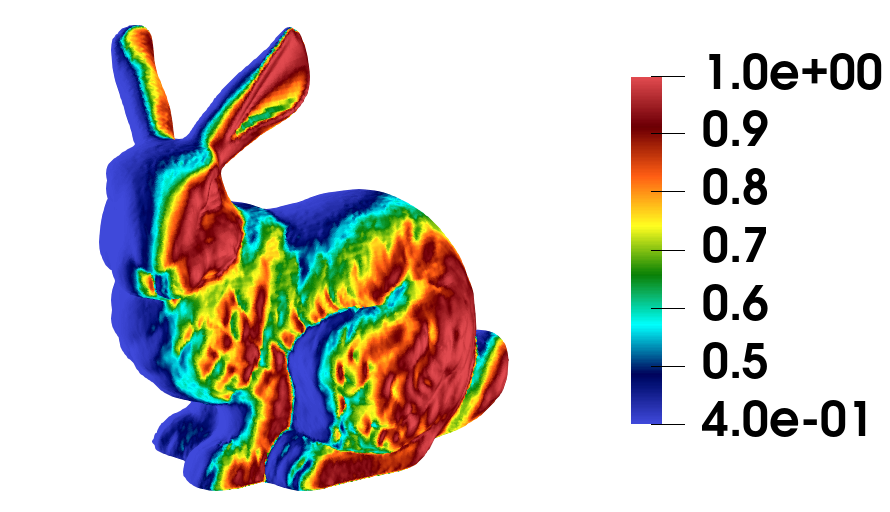}
    \end{subfigure}
    \begin{subfigure}[b]{0.33\textwidth}
        \centering
        \includegraphics[width=0.85\textwidth]{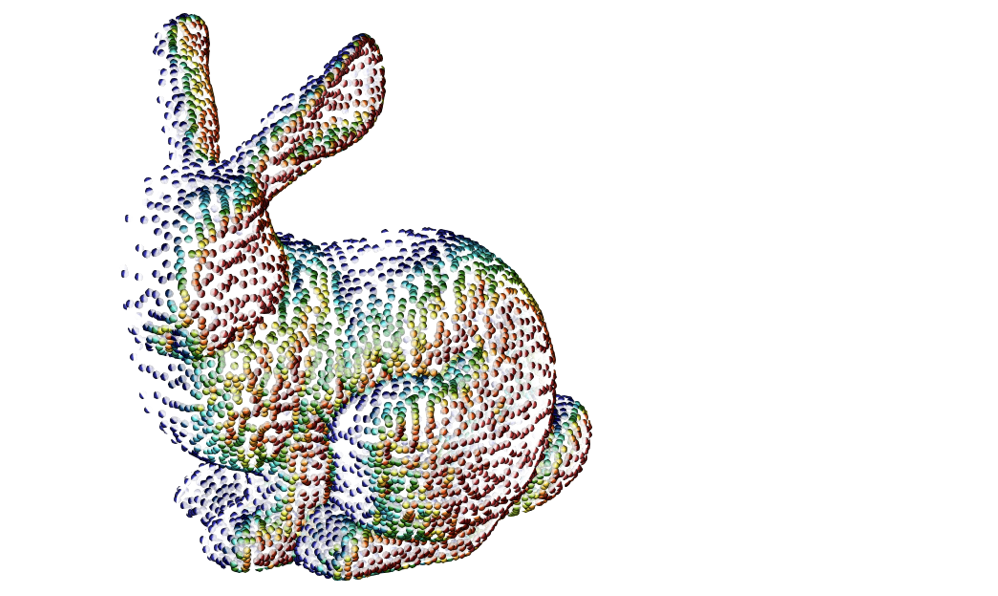}
    \end{subfigure}
    \hspace{-0.5 cm}
    \begin{subfigure}[b]{0.33\textwidth}
        \centering
    \includegraphics[width=0.9\textwidth]{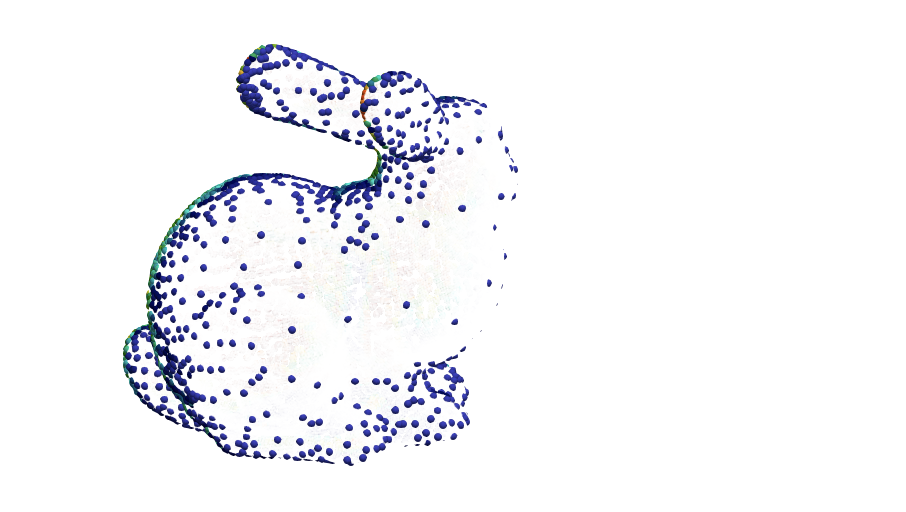}
    \end{subfigure}
    \vskip\baselineskip  
    \begin{subfigure}[b]{0.33\textwidth}
        \centering
        \includegraphics[width=0.85\textwidth]{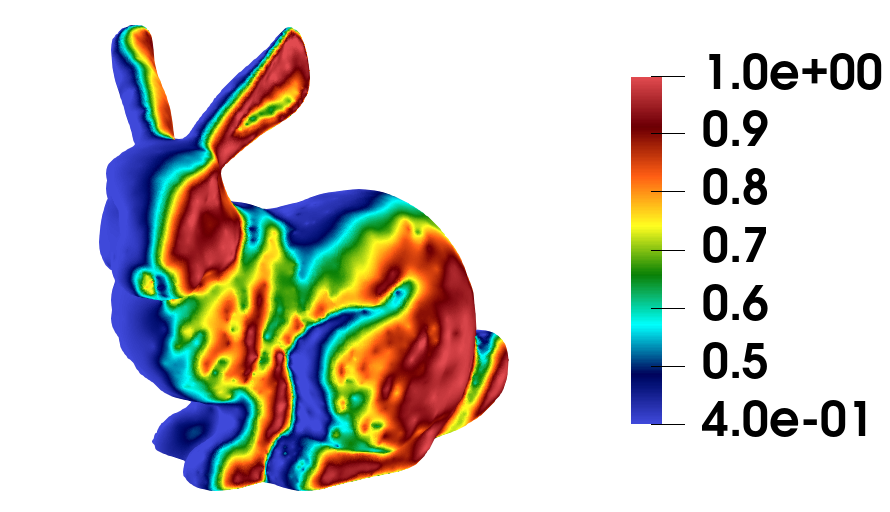}
    \end{subfigure}
    \hspace{-0.2 cm}
    \begin{subfigure}[b]{0.33\textwidth}
        \centering
        \includegraphics[width=0.84\textwidth]{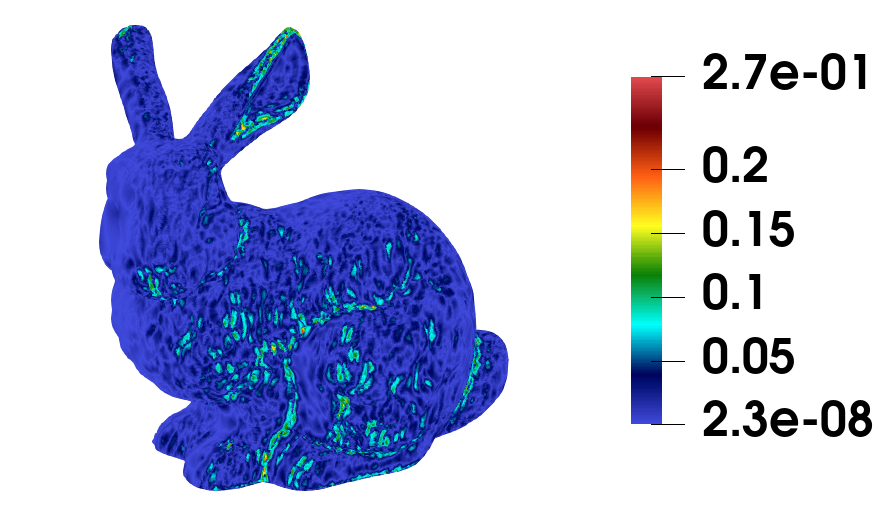}
    \end{subfigure}
    \hspace{-0.2 cm}
    \begin{subfigure}[b]{0.33\textwidth}
        \centering
        \includegraphics[width=0.84\textwidth]{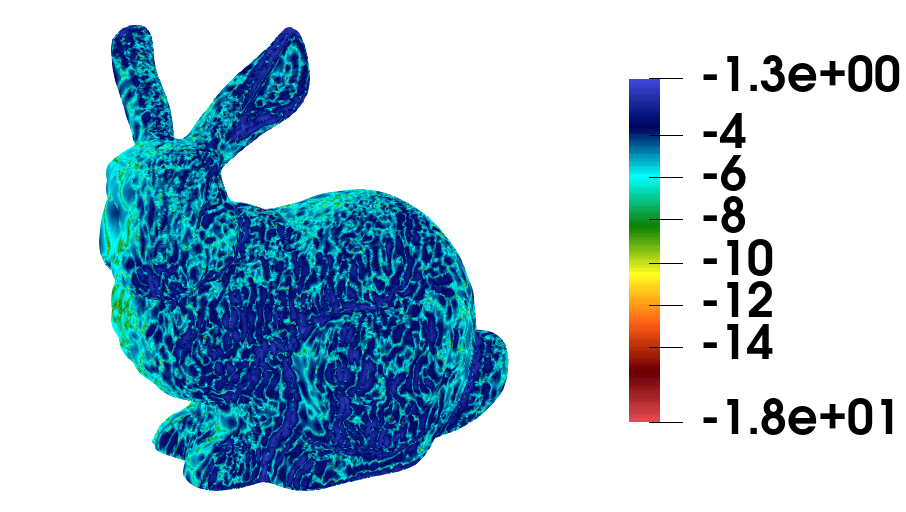}
    \end{subfigure}
\caption{\textbf{Test 4.} Large-scale 3D reflectance on the Stanford Bunny.
First row: reference reflectance field and adaptively selected centers $X_t$ (front and back view).
Second row: reconstructed solution, pointwise absolute error, and pointwise
absolute error in logarithmic {scale}.}
\label{fig:res_test4}
\end{figure}

The selected centers, reconstructed reflectance field, and pointwise errors
are shown in Figure~\ref{fig:res_test4}, while sparsity and reconstruction
accuracy are reported in Table~\ref{tab:test4_results}.
The adaptive sampling concentrates points in regions of high geometric and radiometric variation, enabling an accurate reconstruction from a strongly reduced set of centers.
The proposed framework accurately captures the multiscale reflectance structure demonstrating its scalability
to large-scale three-dimensional problems.

\begin{table}[htbp!]
\centering
\footnotesize
\setlength{\tabcolsep}{12pt}
{\begin{tabular}{@{}lcc@{}}
\toprule
 & $\|\bm\alpha_j\|_0/ |X_t|$ & $e_2$ \\
\midrule
$\bm K_1$ ($\ell_1 = 4.51\cdot10^{-3}$) & $494 / 4857$  & \multirow{3}{*}{$3.77\cdot10^{-2}$} \\[4pt]
$\bm K_2$ ($\ell_2 = 3.35\cdot10^{-2}$) & $1975 / 4857$ & \\[4pt]
$\bm K_3$ ($\ell_3 = 2.49\cdot10^{-1}$) & $86 / 4857$   & \\
\bottomrule
\end{tabular}}
\caption{\textbf{Test 4.} Sparsity of the coefficient vectors $\bm\alpha_1$--$\bm\alpha_3$
across the three kernel lengthscales and relative $\ell^2$ reconstruction error
for the large-scale 3D surface reconstruction ({$|X_t| = 4857$} centers).}
\label{tab:test4_results}
\end{table}

\section{Conclusions}
\label{sct:conclusion}
We have presented a multiresolution framework for sparse
$\ell^1$-regularized kernel approximation that combines samplet-based
compression, native space driven adaptive subsampling and a stabilized
trust-region semismooth Newton solver.
The proposed approach achieves accurate reconstructions while
drastically reducing the number of active coefficients in the solution. The
adaptive selection concentrates data sites in regions of high variability,
leading to reduced problem sizes that reflect the intrinsic multiscale
structure of the data. Furthermore, the samplet compression of the
kernel matrix enables efficient matrix operations, while the TR-SSN method
with the online SVD as interior solver remains stable, 
even in the presence of ill-conditioned kernel matrices.
The method yields significantly sparser representations and improved
computational efficiency compared to single scale formulations and enables
scalability to problem sizes beyond what is feasible with
non-adaptive approaches. In multi-kernel settings, the 
$\ell^1$-regularization acts as an effective mechanism for automatic
kernel selection across different lengthscales.

\section*{Funding}
The authors thank the Swiss National Science Foundation (SNSF)
for the financial support through the grant 215528, ``Large-scale 
kernel methods in financial economics'' and through
the SNSF starting grant ``Multiresolution
methods for unstructured data'' (TMSGI2\_211684).
C. Segala ia a member of the 
Italian National Group of Scientific Calculus (Indam GNCS).

\bibliographystyle{abbrv}
\bibliography{literature}

\end{document}